\documentclass[aip,jmp]{revtex4-1}
\usepackage{accents} 
\usepackage[normalem]{ulem}
\usepackage{cancel}
\usepackage[colorlinks=true]{hyperref}
\hypersetup{urlcolor=blue,citecolor=red}
\usepackage[active]{srcltx}

\newcommand\beq[1]{ \begin{equation}\label{#1} }
\newcommand{\eeq}{ \end{equation} }
\newcommand{\beqno}{ \[ }
\newcommand{\eeqno}{ \] }
\newcommand\beqa[1]{ \begin{eqnarray} \label{#1}}
\newcommand{\eeqa}{ \end{eqnarray} }
\newcommand{\beqano}{ \begin{eqnarray*} }
\newcommand{\eeqano}{ \end{eqnarray*} }
\newcommand\arr[1]{\left\{\begin{array}{l}#1\end{array}\right.}

\newtheorem{theorem}{Theorem}[section]
\newtheorem{definition}{Definition}[section]
\newtheorem{proposition}{Proposition}[section]
\newtheorem{lemma}{Lemma}[section]
\newtheorem{sublemma}{Sublemma}[section]
\newtheorem{remark}{Remark}[section]
\newtheorem{notationalremark}{Notational Remark}[section]
\newtheorem{corollary}{Corollary}[section]
\newtheorem{assumption}{Assumption}[section]
\newtheorem{claim}{Claim}[section]

\newtheorem{tools}{$\negsp\negsp$}[subsection]

\newcommand\thm[1]{ \begin{theorem}\label{#1}}
\newcommand\thmtwo[2]{ \begin{theorem}[#1]\label{#2}}
\newcommand\ethm{ \end{theorem} }
\newcommand\dfn[1]{ \begin{definition}\label{#1} \rm}
\newcommand\dfntwo[2]{ \begin{definition}[#1]\label{#2} \rm}
\newcommand\edfn{ \end{definition} }
\newcommand\pro[1]{ \begin{proposition}\label{#1}}
\newcommand\protwo[2]{ \begin{proposition}[#1]\label{#2}}
\newcommand\epro{ \end{proposition} }
\newcommand\lem[1]{ \begin{lemma}\label{#1}}
\newcommand\lemtwo[2]{ \begin{lemma}[#1]\label{#2}}
\newcommand\elem{ \end{lemma} }
\newcommand\sublem[1]{ \begin{sublemma}\label{#1}}
\newcommand\sublemtwo[2]{ \begin{sublemma}[#1]\label{#2}}
\newcommand\esublem{ \end{sublemma} }
\newcommand\rem[1]{ \begin{remark}\label{#1} \rm}
\newcommand\erem{ \end{remark} }
\newcommand\notrem[1]{ \begin{notationalremark}\label{#1} \rm}
\newcommand\enotrem{ \end{notationalremark} }
\newcommand\cor[1]{ \begin{corollary}\label{#1}}
\newcommand\cortwo[2]{ \begin{corollary}[#1]\label{#2}}
\newcommand\ecor{ \end{corollary} }
\newcommand\asmp[1]{ \begin{assumption}\label{#1}}
\newcommand\asmptwo[2]{ \begin{assumption}[#1]\label{#2}}
\newcommand\easmp{ \end{assumption} }
\newcommand\clm[1]{ \begin{claim}\label{#1}}
\newcommand\eclm{ \end{claim} }

\newcommand\equ[1]{{\rm (\ref{#1})}}

\newcommand{\ie}{{\rm i.e.}}

\newcommand\ovl[1]{ \overline {#1} }

\newcommand{\tr}{{\, {\rm tr}\, }}

\newcommand{\diag}{{ \, \rm diag \, }}

\newcommand{\torus}{ {\mathbf T}   }
\renewcommand{\natural}{ {\mathbf N}   }
\newcommand{\real}{ {\mathbf R}   }
\newcommand{\integer}{ {\mathbf Z}   }
\newcommand{\complex}{ {\mathbf C}   }
\newcommand{\rational}{ {\mathbf Q}  }

\renewcommand{\a }{ {\alpha}   }
\renewcommand{\b}{ {\beta}   }
\newcommand{\g}{ {\gamma}   }

\newcommand{\D}{ {\Delta}   }

\renewcommand{\k}{ {\kappa}   }
\renewcommand{\l}{ {\lambda}   }
\renewcommand{\L}{ {\Lambda}   }
\newcommand{\m}{ {\mu}   }
\newcommand{\n}{ {\nu}   }

\newcommand{\p}{ {\pi}   }
\renewcommand{\P}{ {\Pi}   }

\newcommand{\s}{ {\sigma}   }

\renewcommand{\t}{ {\tau}   }
\newcommand{\f}{ {\varphi}   }

\renewcommand{\o}{ {\omega}   }
\renewcommand{\O}{ {\Omega}   }

\newcommand{\const}{{\, \rm const\, }}

\newcommand{\cA}{ {\cal A} }
\newcommand{\cB}{ {\cal B} }

\newcommand{\cT}{ {\cal T} }
\newcommand{\cR}{ {\cal R} }
\newcommand{\cH}{ {\cal H} }
\newcommand{\cK}{ {\cal K} }
\newcommand{\cC}{ {\cal C} }
\newcommand{\cD}{ {\cal D} }

\newcommand{\cG}{ {\cal G} }
\newcommand{\cL}{ {\cal L} }
\newcommand{\cM}{ {\cal M} }

\newcommand{\cP}{ {\cal P} }
\newcommand{\cI}{ {\cal I} }
\newcommand{{\cJ}}{ {\cal J} }

\newcommand{\cS}{ {\cal S} }

\newcommand\ppu{{ (1) }}
\newcommand\ppd{{ (2) }}
\newcommand\ppt{{ (3) }}

\newcommand\ppj{{ (j) }}

\newcommand\ppi{{ (i) }}

\newcommand\ppo{{ (0) }}

\newcommand\mm{{\rm m}}
\newcommand\MM{{\rm M}}
\usepackage{color}
\definecolor{yellow}{rgb}{0.99, 0.93, 0.0}
\definecolor{applegreen}{rgb}{0.55, 0.71, 0.0}
\definecolor{amber(sae/ece)}{rgb}{1.0, 0.49, 0.0}
\definecolor{amethyst}{rgb}{0.6, 0.4, 0.8}

\newcommand\GG{{\rm G}}

\newcommand\ZZ{{\rm Z}}

\newcommand\RR{{\rm R}}

\newcommand\PP{{\rm P}}
\newcommand\ii{{\rm i}}
\newcommand{\T}{\mathbb{T}}

\newcommand\meas{{\, \rm meas\,}}

\newcommand\HH{{\rm H}}
\newcommand\II{{\rm I}}
\newcommand\FF{{\rm F}}
\newcommand\kk{\textrm{\rm k}}
\newcommand\hh{{\rm h}}
\newcommand\CC{{\rm C}}

\begin{document}

\title{On the co-existence of maximal and whiskered  tori in  the planetary three-body problem} 



\author{Gabriella Pinzari}
\email[]{gabriella.pinzari@math.unipd.it}
\affiliation{Dipartimento di Matematica Tullio Levi--Civita, Universit\`a di Padova, Italy}



\begin{abstract}
In this paper we discuss about the possibility of 
{\it coexistence} of stable and unstable quasi--periodic {\sc kam} tori  in a region of phase space of the three-body problem. 
The  {argument of proof} goes along {{\sc kam} theory and, especially,} the production of two non smoothly related systems of canonical coordinates in the same region of the phase space, the possibility of which is foreseen, for ``properly--degenerate'' systems, by a theorem of Nekhorossev  and Mi{\v{s}}{\v{c}}enko and Fomenko. 
The two coordinate systems are alternative to the classical reduction of the nodes by Jacobi, described, e.g., in~[V.I.~Arnold, Small denominators and problems of stability of motion in classical and celestial mechanics, 18, 85 (1963); p. 141].
\end{abstract} 

\pacs{02. Mathematical methods in physics}

\maketitle 

\section{\label{intro}Introduction}

At the end of the XIX century, motivated by the study of a three--body problem of celestial mechanics, H. Poincar\'e conjectured that a non-integrable system possesses, very close one to the other,  stable and unstable trajectories, ~Ref.\cite[Vol III, Chapter 33, n. 397, p. 389]{poincare}. 
Numerical evidence of such an occurrence has been provided,  since the 50s, as soon as computers could be used to
simulate solutions of differential equations, by   breakthrough papers by Fermi--Pasta--Ulam, H\'enon--Heiles, Izrailev--Chirikov, $\ldots$. The phenomenon was soon understood to be relevant for physics,
since such papers  revealed its appearance
even in regular (e.g. analytic) systems.
{At this respect, the recent development of Aubry--Mather theory allowed to obtain,  under strong assumptions\cite{Note00}, {still failing to be applied to celestial mechanics,} the rigorous proof of existence of motions with {\it any} prescribed frequency.
}
In this paper, we address the question in the case of the celestial  three--body problem; precisely, its {\it planetary} version. 
This is the 4 degrees of freedom problem of three point masses interacting through gravity, where one of the masses (the ``star'') is much larger than the two others (the ``planets'').
We move in the framework of Kolmogorov--Arnold--Moser ({\sc kam}) theory;~Refs.\cite{arnold63c, kolmogorov54, moser1962} and Note\cite{Note0}. {\sc kam} theory has been successfully  applied to problem of celestial mechanics it since the 60s. 
Under the point of view of {\sc kam} theory, the question might be rephrased as wether one can prove coexistence, in a region of phase space,
of quasi--periodic motions, both
  {\it maximal} and {\it whiskered}\cite{Note2} possibly separated, according to Poincar\'e's 
picture and Aubry--Mather theory,  by chaotic regions.
This is precisely the question  {to which this paper is addressed}: we determine {a physical situation}
 where such co--existence
 {is expected to occur}. We consider the following situation, which we shall refer to as   {\it outer retrograde configuration} ({\sc orc}): {\it two planets  describe  almost  co-planar orbits. The outer planet has a retrograde motion\cite{Note1}.}     

{Before describing our result, let us make a digression on the specific features planetary systems and previous literature. }
 It is known that  the {\it two-body problem}, \ie, the problem of the motions
of two  point masses interacting via a law proportional to their inverse squared distance, has, for an open {set}
 of initial conditions, {\it periodic motions} rather than
, more generally, {\it quasi-periodic}. For this system, periodicity consists in the fact that the bounded motions evolve  (according to Kepler's laws) on ellipses,  and are governed by just {\it one} frequency $\n$ proportional to $a^{-3/2}$
where $a$ is the semi-major axis of the ellipse. This pretty remarkable fact unavoidably reflects -- as already underlined by V.~I.~Arnold  in his  1963's paper~Ref.\cite{arnold63} -- on the study of the dynamics of the so-called {\it planetary problem}, i.e., the problem of $(1+N)$ point masses, one of which (``sun'') is of ``order one'', while the remaining $N$ (``planets'') are of much smaller size, interacting through gravity. Indeed, {when the reciprocal attraction among the planets, which is of much smaller order compared to the attraction between any planet and the sun, is set to zero}, the planetary problem reduces to $N$ uncoupled two--body problems (``unperturbed problem'').
He considered the case of $N$ planets in {\it prograde} configuration; i.e., revolving  in the same verse, even though
the question of the sense of rotation, at his time, was definitely of secondary importance, compared to
the 
difficulties that had to be overcome and that we are going to recall.\\
The lack of frequencies{(a translationally invariant  system with $1+N$ bodies possesses $3N$ degrees of freedom. In the case of the planetary problem, it exhibits, as mentioned, only $N<3N$ frequencies)}  in the unperturbed problem  was named by Arnold  {\it proper degeneracy}. It represented  a serious difficulty, if one wanted (as he was aiming to do) to apply  Kolmogorov's theorem,~Ref.\cite{kolmogorov54} to the planetary problem.   
\\
At a technical level, the  appearance of the proper  degeneracy consists, we might say, of a  ``loss of frequencies'', caused by the  ``too many'' (or, better {\it Poisson non commuting} \cite{Note3}, see below) first integrals of motion. For such abundance, this kind of systems is often called {\it super--integrable}. Despite of the fact that 
the solutions of the two-body problem are known since Newton's times,  a general, theoretical setting clearly explaining the phenomenon  has been given only recently, thanks to the works by Nekhorossev  and Mi{\v{s}}{\v{c}}enko and Fomenko,~Refs.\cite{nehorosev72, MishchenkoF81} (hereafter, {\sc nmf}). 
The three authors proved a generalization of the best known Liouville-Arnold theorem,~Ref.\cite{arnold63a}
 which clearly relates  the loss of frequencies to the existence of Poisson non commuting independent integrals. They proved that, to an integrable Hamiltonian system with $n$ degrees of freedom which, in addition to  $n$ independent and commuting first integrals, affords  additional, {independent, from the first $n$,} $n_1\le n$  integrals which do not commute with {\it all}
 the integrals of the first family, one can associate canonical coordinates including only $n_0:=n-n_1$
  action--angle {pairs} $(\II,\f)=(\II_1,\cdots, \II_{n_0},\f_1,\cdots,\f_{n_0})$ (analogous to the ones of Arnold--Liouville case), and, in addition, certain other couples $z=(p,q)=(p_1,\cdots, p_{n_1}, q_1,\cdots, q_{n_1})$, usually referred to as {\it degenerate} coordinates. The degenerate coordinates  {\it are not uniquely defined}, and this is precisely the aspect  
that, in this paper, we shall exploit. \\
Indeed, a dynamical system that is close to a super-integrable system may be written as
\beq{proper degeneracy}\HH(\II,\f,p,q)=\hh(\II)+\m f(\II,\f,p,q)\eeq
where $(\II,\f,p,q)$ is one of the various (as foreseen by {\sc nmf} Theorem) sets of canonical coordinates associated with the unperturbed super-integrable term $\hh$. Now, while, given the $\II$'s, $\hh$ is uniquely determined, the form of $f$, instead, strongly depends on the choice of coordinates. On the other hand,  it is known since Arnold's paper~Ref.\cite{arnold63} that, for system of the form~\equ{proper degeneracy},  $f$ may have a strategic importance.

As an outstanding example, let us recall just the case considered by Arnold in~Ref.\cite{arnold63}. He wanted to prove (via an application of Kolmogorov's theorem) the existence of plenty of quasi--periodic, maximal tori, forming
a positive measure set in phase space. He announced the result (known as ``Arnold Theorem'') at the 1962 ICM.
Clearly, such result was
 going in the direction of the proof of stability of the Solar System, and for this Kolmogorov and Arnold were awarded, in 1965, of the Lenin prize. 
 However, in order to obtain such result he was aware that he had to overcome the problem of
the lack of frequencies in the unperturbed part (indeed succeeding in this), but this was not the only one. As for the choice of coordinates, Arnold considered, in the case of the planar problem, {\it Poincar\'e coordinates}, as described in~Ref.\cite[Chapter III, \S 2, n.4]{arnold63}.  In term of such coordinates, the 
Hamiltonian of the planetary problem takes the form in~\equ{proper degeneracy}, with
 $n_0$ equal to  the number of planets $N$, $n_1=N$ (so that the total number of degrees of freedom in $\real^2$ is $2N$),  $(\II$, $\f):=(\L,\l)\subset \real^N\times \torus^N$ (where $\torus:=\real/(2\p\integer)$) suitable action--angle couples related to the semi--major axis and the area spanned by the ellipse, $z=(p,q):=(\eta,\xi)\subset \real^{N}\times\real^{N}$ suitable degenerate coordinates related to  the orientation of such ellipses, $\hh=\hh_{\rm k}$ the  Keplerian Hamiltonian; $\m$ a small a--dimensional parameter measuring the maximum planet/star mass ratio and, finally, $f(\II,\f,p,q)=f_{\rm Poin}(\L,\l,\eta,\xi)$ a perturbing function related to the small mutual interactions among planets. Arnold observed that  the average value $\ovl {f_{\rm Poin}}(\L,\eta,\xi)$ with respect to the $\l$'s of the perturbing function $f_{\rm Poin}(\L,\l,\eta,\xi)$
 by symmetry reasons, has an {\it elliptic equilibrium point} for $(\eta,\xi)=0$ (corresponding to circular motions of the planets around their sun), for all $\L$. So he managed to construct, for degenerate systems   of the form~\equ{proper degeneracy} with $\ovl f(I,z)$  having an elliptic equilibrium in $z=0$ for all $\II$, a careful version of Kolmogorov Theorem,~Ref.\cite[Fundamental Theorem]{arnold63}    {based on a generalized non--degeneracy condition (``full torsion''), inspired to Kolmogorov,  according to which one should check, besides of the non--singularity of the Hessian matrix {$\partial_{\II_i\II_j}^2\hh(I)$ in  \equ{proper degeneracy}}
 also the one of the matrix of the coefficients of the second--order term  of the Birkhoff normal form associated to the elliptic equilibrium (see~Ref.\cite{hoferZ94})}.
Arnold successfully applied his Fundamental Theorem to the case of the planar problem with $N=2$ planets. However, while the extension to the planar problem with a generic number of planets revealed to be straightforward,~Ref.\cite{pinzari-th09} (see~Ref.\cite{fejoz04} for a previous result with a different strategy), the treatment of the  problem in space contained strong extra-difficulties.
Indeed, switching from planar to spatial Poincar\'e coordinates, the  averaged perturbing function $\ovl {f_{\rm Poin}}$ still exhibits an elliptic equilibrium in correspondence {\it circular and co--planar} motions, but such equilibrium is {\it degenerate}, in the sense that the eigenvalues of the quadratic part of  $\ovl {f_{\rm Poin}}$ verify, identically, two linear combinations with integer coefficients (known in the field as {\it secular resonances}). A fact strongly preventing, in principle, {the construction of the Birkhoff normal form and hence} the possibility of checking {the full torsion condition.}  But this is not all: 
{a formal evaluation of the torsion, attempted asymptotically by M. Herman,~Ref.\cite{herman09},
seemed to  suggest, in absence of proper reductions of the rotation invariance, an identically vanishing determinant (implying the impossibility of applying the Fundamental Theorem to the general problem), a fact next rigorously proved, by L. Chierchia and the author, in~Ref.\cite{chierchiaPi11c}}.\\
It may be argued that Arnold  felt that a difficulty of this kind might appear, since, without explaining his motivations, in~Ref.\cite[p. 141--42]{arnold63}, he  suggested to ``change coordinates'', without going further.
Completion of the proof of his theorem revealed it to be  more difficult than expected, and the story reached a conclusion only fifty years later, thanks to contributions by J. Laskar, P. Robutel, M. Herman, J.F\'ejoz,  L. Chierchia and the author,~Refs.\cite{laskarR95, robutel95,fejoz04,pinzari-th09, chierchiaPi11b}. Comprehensive reviews appeared in~Refs.\cite{fejoz13,chierchiaPi14}, to which papers we refer  the interested reader.
For the purposes of this paper we only mention that the solution Arnold had in mind, based on changing coordinates was considered, formally, in a particular case, by Malige, Robutel and Laskar,~Ref.\cite{maligeRL02}, and next completely achieved by the author,~Ref.\cite{pinzari-th09}, published in~Refs.\cite{chierchiaPi11b,pinzari-th09}.
The  long proof of Arnold Theorem  should give, we hope,  an idea that, from a practical point of view, 
producing ``good'' canonical coordinates, which should: (i)  leave the unperturbed part unvaried; (ii)
 overcome the degeneracies caused by SO(3) invariance and, eventually, (iii)  preserve symmetries,  parities, equilibria $\ldots$ from which to depart in order to apply a perturbative scheme (e.g., in the case of Arnold Theorem, the Fundamental Theorem developed around  the elliptic equilibrium), is other than ``easy'' or ``straightforward''. \\
 In this paper, we  {look at } the three-body problem in the {\sc orc} configuration {by means of (basically)}  {\it two} sets of canonical coordinates.
The former of such two systems of coordinates is a modification of the {so--called  ``regular, planetary and symplectic'' ({\rm rps}) coordinates}, proposed in~Refs.\cite{pinzari-th09, chierchiaPi11b}. The latter, called ``perihelia reduction'' ({\rm p}) has been proposed in~Ref.\cite{pinzari15}. {Both such systems of coordinates describe regularly }
 co--planar motions{,  which evolve on suitable invariant manifolds of each phase space}.
{Each of} such invariant manifolds  {turns to be an} equilibrium {for suitable  truncated and averaged systems (where the average is performed with respect to fast angle coordinates in each of such sets),} parametrized by the value of {certain other} action coordinates, {which play the r\^ole of quasi--integrals of motion}.\\
We provide the complete proof of the existence of a positive measure set with a maximal number of frequencies for the full system, both in the case of planar and spatial problem. More precisely, we  prove the following result (a more precise formulation  will be given in course of the paper; see  Theorem~\ref{stable tori}).
\vskip.1in
{\bf  Theorem} {\it There exists an eight-dimensional open region 
of  phase space{,} contained in the holomorphy domain of the Hamiltonian{,}  almost completely filled with a positive measure set of quasi--periodic motions with {maximal number of} incommensurate frequencies. The motions on such tori are in {\sc orc}.}

\vskip.1in
{
The proof of the Theorem adapts the techniques of~Refs.\cite{chierchiaPi11b, pinzari-th09}, which, as recalled above, dealt with the {\it prograde} case. {Although} the strategy is the same, nevertheless,  certain structural differences between the two settings
do exist, which we  point out. The most remarkable one is related to the effect of the rotation invariance 
in the two cases, in relation with the elliptic nature of the co--planar, co--circular equilibrium. While, in the case treated in~Refs.\cite{chierchiaPi11b, pinzari-th09}, such ellipticity is a mere consequence of its invariance by reflections and rotations (a fact already known to M. Herman,~Ref.\cite{herman09}), in the {\it retrograde} case, it is not so, but needs to be checked specifically.}

 {We conclude this introduction with recalling related literature.}\\
 {The existence of a positive measure set of Lagrangian tori with maximal number (see the next section) of frequencies for the general planetary problem, in the regime of well spaced orbits, small eccentricities and small inclinations, has been established in the papers~Refs.\cite{arnold63, laskarR95, fejoz04, pinzari-th09, chierchiaPi11b, pinzari13,  pinzari15}. We refer to such technical papers for details, to~Refs.\cite{fejoz13, chierchiaPi14} for reviews.
However, the cases treated in the literature above, even though containing all the necessary information, are not perfectly suited to the proof of the the {\bf Theorem}.}

{The papers in~Refs.\cite{arnold63,  fejoz04, pinzari-th09, chierchiaPi11b, laskarR95, robutel95} deal with maximal quasi--periodic tori in the case when the planets revolve all in the same verse, and eccentricities and inclinations are small.
 The invariant set (so-called ``{\it Kolmogorov-set}\,'') {for the (rescaled) Hamiltonian in~\equ{helio}} is proved to fill almost completely the  {region in phase space with small eccentricities and inclinations}, up to a residual set with measure going to zero with {the parameter $\m$ in~\equ{helio} and with} the maximum {$\varepsilon$} of eccentricities and inclinations. 
A (maybe optimal) estimate about the strength at which  such measure goes to zero is provided in~Refs.\cite{pinzari-th09, chierchiaPi11b}. }

{In~Refs.\cite{pinzari13, pinzari15} maximal quasi--periodic tori have been constructed out of the small eccentricities and inclination constraint. In such papers, the measure of the Kolmogorov-set has been found to increase while the masses  decrease and the mutual semi--major axes ratios increase, independently of the values of eccentricities and inclinations. A suitable constraint on the semi--axes ratios is however imposed. Finally, the sense of revolution of the planets is the same for all of them. }

{A  study of quasi--periodic motions (including retrograde ones)  bifurcating  from relative equilibria 
   appeared in~Refs.\cite{palacianSY13, palacianSY15}. {However, the measure estimates obtained in those papers, based on a bit different framework (Withney regularity and no use of Birkhoff normal form; see, e.g., Ref.\cite[Theorem 5.1]{palacianSY15}),   are not suited to the purposes of the paper, where the point of view is closer to Ref.\cite{arnold63}.}}

\section{\label{Existence of maximal tori} Set Up}

{The three--body problem is the   dynamical system formed by three point masses in $\real^3$ interacting through gravity only. 
The system has thus nine degrees of freedom, meaning with this that its evolution is described by a system of differential equations  having order eighteen. However, it is also well known that the system possesses several constant motions, and, even though the number of such constant is not sufficient (as Poincar\'e showed, Ref.\cite{poincare}) to guarantee integrability by quadratures, 
nevertheless, it allows to reduce the order
of the equations from eighteen to eight. This complete reduction was firstly considered by
Jacobi and refined by R. Radau, Refs.\cite{jacobi1842, radau1868} (see the next section). The first step to achieve it consists of  getting free of the translation invariance, caused by the conservation of the total linear momentum (the velocity of the center of mass of the system). In literature one finds essentially two ways to do it, usually referred to as ``Jacobi'' or ``heliocentric'' coordinates. They both can be described as  linear changes
of coordinates, if the Hamiltonian of the system is initially written in impulse--position coordinates. We refer to Refs.\cite{arnold63, robutel95} for a complete description of them.
 According to  the heliocentric reduction, in a system
where the masses are denoted as $m_0$,  $\m\,m_1$, $\m\,m_2$, where $\m$ is a prefixed pure number  (the case  $\m\ll 1$ and $m_i$'s of the same strength being usually referred as ``planetary'' problem), the motions  are described the Hamilton equations of the six degrees of freedom Hamiltonian
}
\beq{helio}{\rm H}{(y,x)}=\frac{|y^\ppu|^2}{2{\rm  m}_1}-\frac{{\rm  m}_1{\rm  M}_1}{|x^\ppu|}+\frac{|y^\ppd|^2}{2{\rm  m}_2}-\frac{{\rm  m}_2{\rm  M}_2}{|x^\ppd|}+\m\Big(-\frac{m_1m_2}{|x^\ppu-x^\ppd|}+\frac{y^\ppu\cdot y^\ppd}{m_0}\Big)\eeq
 where 
 \beq{reduced masses}
{\rm  m}_i:=\frac{m_0 m_i}{m_0+\m m_i}=m_i+{\rm O}(\m)\qquad {\rm  M}_i:=m_0+\m m_i=m_0+{\rm O}(\m)
\eeq
    are the ``reduced masses''; $y^\ppi\in \real^3$, $x^\ppi\in \real^3$, and the collision set
    $$\D:=\Big\{x^\ppu=0, {\rm or}\ x^\ppd=0, {\rm or}\ x^\ppu=x^\ppd\Big\}$$
    is to be excluded. {Incidentally, the two terms in the perturbing function are sometimes referred  to as {\it direct (or Newtonian), indirect part}, respectively.
}

\vskip.1in
    {After the linear momentum reduction, the next issue is to get rid of rotation invariance of the Hamiltonian \equ{helio}, caused by the conservation of} the three components $\CC_1$, $\CC_2$ and $\CC_3$ of total angular momentum of the system:
\beqa{tot ang mom}\CC=\CC^\ppu+\CC^\ppd\qquad{\rm with}\qquad \CC^\ppi:=x^\ppi\times y^\ppi\ .\eeqa
{This further step is more subtle than the previous one for two reasons. The first obvious reason is that, differently from the linear momentum reduction, it cannot be obtained via a linear transformation. But the main reason is that, in the case of the problem in $\real^3$,  the $\CC_i$'s do not Poisson--commute, one cannot think to a elimination ``by quadratures'' (i.e., one cannot think, roughly speaking, to use them as generalized momenta, mimicking the linear momentum reduction procedure).  The only widely known method in the field (recalled in the next section) is due to C. G. Jacobi, R. Radau and, after the work of A. Deprit (see the next section), it is available for any number of particles.  Another reduction, called ``reduction of perihelia'', has been recently proposed in Ref.\cite{pinzari15} and, as well as the previous one, is available for any number of bodies. It will be recalled  in Section~\ref{generalities} in the particular case of three bodies.  For the results of the paper, both the mentioned reductions will be used. A new unified proof of their canonical character is presented in Section~\ref{appendix D}.}

\subsection{\label{The canonical setting}The Jacobi--Radau--Deprit coordinates}
In the XIX Century C.G.Jacobi,~Ref.\cite{jacobi1842}, found a tricky procedure, {that is usually referred to as ``Jacobi's reduction of the nodes''}  that allowed him to write the
differential equations of the spatial three--body problem {as a system of order eight.} {His speculations were refined by R. Radau, Ref.\cite{radau1868}, who wrote such equations as a system of eight equations of order one, corresponding to the Hamilton equations of a four degrees of freedom Hamiltonian. Even though the original Jacobi--Radau's work was suited for a general two--particles system enjoying rotation invariance,  it is customary (compare, e.g.,  Ref.\cite{arnold63}) to refer with the same  name a slightly modified procedure such in a way that the  --integrable -- (translationally reduced) two--body
terms 
 \beq{K***}\frac{|y^\ppj|^2}{2{\rm  m}_j}-\frac{{\rm  m}_j{\rm  M}_j}{|x^\ppj|}\eeq
appearing in the Hamiltonian \equ{helio}, are put in action--angle form, in the sense of Liouville--Arnold Theorem, Ref.\cite{arnold63a}.
As known, the Liouville--Arnold form for \equ{K***} is one--dimensional
\beq{h(L)}{\rm h}_{{\rm  k}}^\ppj(\L_j):=-\frac{{\rm  m}_j^3{\rm  M}_j^2}{2\L_j^2},\eeq
with the action $\L_j$'s being related to the semi--major axis $a_i$ of the Keplerian ellipse via 
\beq{aL}
\L_j={\rm m}_j\sqrt{{\rm M}_j a_j}\ .
\eeq
Jacobi's trick consisted in fixing in advance a rotating reference frame
having its third axis in the (constant) direction of the total angular momentum ${\rm C}$ and its first axis in the (moving) direction of the so--called ``nodes lines''. Astronomers call so the straight line determined by the intersection (provided it is well defined) of the instantaneous planes of of the orbits of the two planets; i.e., the planes $\P_j(t)=(y^\ppj(t), x^\ppj(t))$, $j=1$, $2$. Jacobi and Radau proved that, even though the reference frame moves, nevertheless, the system of eight coordinates 
given by the  ``planar Delaunay elements''  (see Ref.\cite{GG83})
\beqa{j***}{\rm j}=(\L_1,\L_2,\GG_2,\GG_2, \ell_1, \ell_2, \g_1, \g_2)\eeqa
on the moving planes $\P_1(t)$, $\P_2(t)$
induces an injection
\beqa{j}\phi_{\rm j}:\quad {\rm j}=(\L_1,\L_2,\GG_2,\GG_2, \ell_1, \ell_2, \g_1, \g_2)\in \real^4\times \torus^4\to \cC=(y^\ppu, y^\ppd, x^\ppu, x^\ppd)\in \real^{12}\eeqa
such that the motion of the system are the solutions of the Hamilton Equation of  $\HH_{\rm j}:=\HH\circ\phi_{\rm j}$, which, moreover, turns to depend parametrically only of the length $\GG=|\CC|$ (namely, it does not depend on the direction of $\CC$; a fact vaguely attributed, at that time, to the rotation invariance of the system, but next fully understood thanks to the work of A. Deprit, much years later; see below).
}
The reduction by Jacobi and Radau was well known to V.I.Arnold, who mentioned it in \cite[II, \S5, n.4, p.141]{arnold63} as a unavoidable tool (due to certain degeneracies that appeared if such tool was not used) in order to prove the stability of the planetary, spatial three--body problem in the perturbative setting. {The proof of stability  that he had discovered
consisted of checking the non--singularity of  a certain matrix (``torsion''; see Section~\ref{domain of rps+}) related to a certain averaged perturbing function.
However, having treated with much detail
the case of the planar problem (which had required a considerable amount of computation to him), in order to simplify the analysis, instead of computing the torsion of the  spatial problem directly, he
preferred to try  to reduce the spatial problem to a perturbed planar one,  so as to use the computations he had already done. 
This led him in error (roughly, because Jacobi's reduction is singular for co--planar motions). The computation of the torsion for the spatial problem, using Jacobi's reduction of the nodes, was next completed
 by Laskar and Robutel};~Refs.\cite{laskarR95, robutel95}.\\
How to obtain a generalization of Jacobi--Radau coordinates to the case of more than three bodies has been a one century long problem of mechanics, that Arnold invoked as an obstruction to the extension of his outline of proof to the general planetary problem in the space,~Ref.\cite[Ch. III, \S5, n.5, p. 141]{arnold63}.
A mild progress was offered only in 1982 by F. Boigey who, during her PhD, obtained
a Jacobi--like  reduction for the problem of four bodies,~Ref.\cite{boigey82}. One year after, in 1983, and twenty years after Arnold's work, A. Deprit, positively impressed (as he declared in the introduction) by Boigey's work, discovered a set of canonical coordinates  
that {(in a sense; see below)} reduces to Jacobi--Radau's for $N=2$ and to Boigey's for $N=3$. For some strange reason, Deprit's work was  overlooked\cite{Note3bis} for further twenty years
and Deprit himself seemed to be not much confident about the utility of his coordinates\cite{Note4}.
 During such time, the complete proof of Arnold's Theorem had been obtained by Herman--F\'ejoz,~Refs.\cite{fejoz04, herman09}, with a different strategy (the problem of degeneracies in the Hamiltonian was overcome with an abstract argument of Lagrangian mechanics, without using explicit coordinates). Deprit's work was rediscovered by the author during her PhD, precisely in the framework of obtaining a proof of Arnold's theorem accordingly to the original strategy. The production, in \cite{pinzari-th09} (published in~Ref.\cite{chierchiaPi11b}) of a new set ``regularizing'' coordinates (see also the next section) much similar to Poincar\'e's coordinates, but better suited to bypass the problem of degeneracies, was the final key to reach the objective.

{As in the case of Jacobi's reduction, it  is  customary to call ``Deprit's coordinates'' a modified version  (in fact, the form they were rediscovered in Ref.\cite{pinzari-th09}) of the original set discussed in Ref.\cite{deprit83}, such in a way to satisfy \equ{h(L)}. We recall such modified version, in the case of a system of two particles.}\\
We fix a domain $\cD_{{\rm  jrd}}\subset \real^{12}$ in phase space as follows. Let $(k^\ppu, k^\ppd, k^\ppt)$ be a prefixed orthonormal frame in $\real^3$.  For {the Cartesian coordinates} ${\cC_{art}=}(y^\ppu, y^\ppd,x^\ppu, x^\ppd)\in \cD_{{\rm  jrd}}$, we assume  that
the orbits $t\to (x^\ppj(t), y^\ppj(t))$ generated by the  Hamiltonians~\equ{K***} with initial datum $(x^\ppj, y^\ppj)$ are  ellipses with
 non-vanishing eccentricity, {belonging to different planes, never coinciding with the $(1,2)$ plane}.  
 Then we denote as
 ${\rm P}^{(j)}$  the unit vectors pointing in the directions of the perihelia; as $a_j$ the semi--major axes; 
 {as $\ell_j$ the ``mean anomaly'' of $x^\ppj$(which, we recall, is defined as area of the elliptic sector from ${\rm P}^{(j)}$ to $x^\ppj$ ``normalized at $2\p$'');}
 as  ${\rm C}^\ppj=x^\ppj\times y^\ppj$, $j=1$, $2$,  the angular momenta of the two planets and ${\rm C}:={\rm C}^\ppu+{\rm C}^\ppd$ the total angular momentum integral. {By assumption,} the ``nodes''
\beqa{Dep nodes}\n_1:=k^\ppt\times{\rm C}\ ,\quad \n:={\rm C}\times {\rm C}^\ppu=\CC^\ppd\times\CC^\ppu\eeqa
do not vanish. 
For three vectors $u$, $v$, $w$ with $u$, $v\perp$ $w$, we denote as $\a_{w}(u,v)$  the angle formed by $u$ to $v$ relatively to the positive (counterclockwise) orientation established by $w$.
The  
${\rm jrd}:=\phi_{\rm jrd}^{-1}\big(\cC_{art}\big)$
coordinates
 are defined
via the following formulae: 
\beqa{coordinates}
\begin{array}{llllrrr}
 \left\{
\begin{array}{lrrr}
 Z:={\rm C}\cdot k^\ppt \\
{\rm G}:=\|{\rm C}\| \\
{\rm G}_1:=\|{\rm C}^\ppu\|\\
{\rm G}_2:=\|{\rm C}^\ppd\|\\
\L_j:={\rm  M}_j\sqrt{{\rm  m}_j a_j}
\end{array}
\right.\qquad
\left\{
\begin{array}{lrrr}
\zeta:=\a_{k^\ppt}(k^\ppu, \n_1)\qquad& \\
 \g:=\a_{{\rm C}}(\n_1, \n)\qquad&\\
 {\g}_1:=\a_{{{\rm C}^\ppu}}(\n, {\rm P}^\ppu)&\\
 {\g}_2:=\a_{{\rm C}^\ppd}(\n,{\rm P}^\ppd)&\\
\ell_j := {\rm mean\ anomaly\ of}\ x^\ppj& 
\end{array}
\right.
\end{array}
\eeqa
Note that such definitions, together with  assumptions of non--vanishing of the nodes~\equ{Dep nodes}, imply 
that
\beqa{Dep conditions}
\min\Big\{\GG_1,\ \GG_2,\ \GG\Big\}>0\ ,\quad -\GG<\ZZ<\GG\ ,\quad |\GG_1-\GG_2|<\GG<\GG_1+\GG_2\ .
\eeqa
{
{The coordinates \equ{coordinates} will be denoted in this paper as
\beq{J}\textrm{\rm jrd}:=(\L_1,\L_2,\GG_1,\GG_2,\GG,\ZZ,\ell_1,\ell_2,\g_1,\g_2,\g, \zeta)\eeq
 by the names of their discoverers. Observe, for comparison, that, while \equ{j}
 is a  canonical injection, the \equ{J}'s come out from a canonical change of coordinates 
 \beqa{jrd map}\phi_{\rm jrd}:\quad {\rm jrd}=(\L_1,\L_2,\GG_1,\GG_2,\GG,\ZZ,\ell_1,\ell_2,\g_1,\g_2,\g, \zeta)\to \cC_{\rm jrd}=(y_{\rm jrd}^\ppu,y_{\rm jrd}^\ppd, x_{\rm jrd}^\ppu, x_{\rm jrd}^\ppd)\eeqa
i.e., 
which preserves the standard 2--form
\beqa{jrd canonical} \sum_{j=1}^2(d\L_j\wedge d\ell_j+d\GG_j\wedge d\g_j)+
d\GG\wedge d\g+d\ZZ\wedge d\zeta=\sum_{j=1}^2\sum_{i=1}^3dy^\ppj_{{\rm jrd},i}\wedge dx^\ppj_{{\rm jrd},i}\ .\eeqa
From the practical point of view, there is no difference with the two, because  indeed the coordinates $\ZZ$, $\zeta$ and $\g$, which do not appear in \equ{j}, are cyclic to the Hamiltonian
\beqa{jrd Ham}
\HH_{\rm jrd}:=\HH\circ\phi_{\rm jrd}=\hh_{\rm k}+\m f_{\rm jrd}
\eeqa
(where $\HH$ is as in~\equ{helio})
 which in fact coincides with $\HH_{\rm j}:=\HH\circ\phi_{\rm j}$. 
 {Namely, in both cases, one obtains a 
 a 10--dimensional Hamiltonian that depends on the eight coordinates
in \equ{j***} and G but does not depend on the  DepritÕs angle $\gamma$
conjugated to G} (so that one can regard $\GG$ as a parameter and regard the system as having four degrees of freedom).}
On the other side (and we consider this the major advantage of the  \equ{coordinates}'s),   completing the \equ{j***}'s  with
the quadruplet $(\ZZ, \GG, \zeta, \g)$  is just what is needed to obtain a generalization to any $N\ge 2$ by induction, as shown in Ref.\cite{pinzari-th09} (the original proof by Deprit in Ref.\cite{deprit83} is not by induction). }\\
 {The analytical form of the map \equ{jrd map} has beed provided in Refs.\cite{pinzari-th09, chierchiaPi11b}, in a different 
 framework compared to the original one proposed by Deprit (who, in
  Ref.\cite{deprit83} proposed to use ``quaternions''). See also Appendix~\ref{expression of jrd}.}

\section{Full dimensional quasi--periodic motions in the retrograde problem}
\subsection{Analytical tools}
In this section we provide the analytical tools in order to prove the existence of full dimensional quasi--periodic motions in the retrograde problem. 
{The proof goes through a change of coordinates that transforms the Hamiltonian \equ{helio} to a new Hamiltonian with the average with respect to the angles conjugated to the $\L_i's$ being Birkhoff normal form.
We remark, at this respect that, 
instead of constructing {\it directly} a system of canonical coordinates that do the job and, at the same time, reduce the number of degrees of freedom to {\it four} (accordingly to the previous section),  with the purpose of limiting the computations at a minimum,
we  mimic the procedure followed in Ref.\cite[Sections 4--6]{chierchiaPi11b}, where
the construction of the Birkhoff normal form for the fully reduced system was obtained via the passage to an intermediate system
(``partial reduction'') where the number of effective degrees of freedom is one over the minimum. In particular, we shall establish
a correspondence between the formulae of the ``prograde'' case, treated in  Refs.\cite{pinzari-th09, chierchiaPi11b}, and the ones of the retrograde one, considered here. As well as in the previous literature, the coordinates that we propose have nice ``regularizing properties'', described in the next section.
}

\subsubsection{\label{prograde case} Regularizing coordinates}
The  {map \equ{jrd map}}, and hence the Hamiltonian \equ{jrd Ham}, are not well defined when  some of the following equalities is verified
\beqa{singularities}&&\GG_1=\L_1\quad {\rm or}\quad \GG_2=\L_2\quad {\rm or}\quad\GG=\GG_1+\GG_2\quad {\rm or}\quad\GG=\GG_1-\GG_2\quad {\rm or}\quad\GG=\GG_2-\GG_1\nonumber\\
&&\text{or}\quad \ZZ=\GG\quad\text{or}\quad \ZZ=-\GG\quad {\rm or}\quad\GG_i=0\quad {\rm or}\quad\GG=0
\eeqa
Some advantage is taken  in using regularized coordinates, at expenses of introducing an extra--integral, as now we describe.{
%
 %
 The coordinates doing the job in the case of the prograde problem have been found in Ref.\cite{pinzari-th09} and have been called, in Ref.\cite{chierchiaPi11b}, ``regular, planetary and symplectic'' (rps). To avoid confusions with the coordinates that we shall introduce in a while, the rps coordinates in the case of the three--body problem will be denoted, in the paper, with the symbols
\beqa{rps complex}
 \textrm{\rm rps}&:=&(\L,\widetilde\l,  \widetilde \eta, \widetilde \xi, \widetilde p, \widetilde q, \widetilde P,\widetilde Q)=(\L_1,\L_2,\widetilde\l_1,\widetilde\l_2,  \widetilde \eta_1,\widetilde\eta_2, \widetilde \xi_1,\widetilde\xi_2, \widetilde p,\widetilde q, \widetilde P, \widetilde Q)\eeqa
(they correspond to the coordinates for $n=2$ in Ref.\cite{chierchiaPi11b}, and thereby named  $L_1$, $L_2$, $\l_1$,  $\l_2$, $ \eta_1$, $ \eta_2$, $ \xi_1$, $\xi_2$, $ p_1$, $ q_1$, $ p_2$, $ q_2$, respectively).
{For later convenience, we introduce the complexified version of the \equ{rps complex}'s, defined as}
\beqano
 \textrm{\rm rps}^{\complex}&:=&(\L,\widetilde\l,  \widetilde t, \widetilde t^*, \widetilde T, \widetilde T^*)=(\L_1,\L_2,\widetilde\l_1,\widetilde\l_2,  \widetilde \t_1,\widetilde\t_2, \widetilde \t_3,\widetilde t_1^*, \widetilde t_2^*,\widetilde t_3^*, \widetilde T, \widetilde T^*)\eeqano
where
 \beqa{PROld+}
\widetilde t_1&:=&
  \frac{\widetilde \eta_1-{\rm i}\widetilde \xi_1}{\sqrt2}\qquad \widetilde t_2:=
  \frac{\widetilde \eta_2-{\rm i}\widetilde \xi_2}{\sqrt2}\qquad\  \widetilde t_3:=
   \frac{\widetilde p-{\rm i}\widetilde q}{\sqrt2}\qquad\widetilde T:=
  \frac{\widetilde P-{\rm i}\widetilde Q}{\sqrt2}\nonumber\\\nonumber\\
  \widetilde t^*_1&:=&
  \frac{\widetilde \eta_1+{\rm i}\widetilde \xi_1}{\sqrt2{\rm i}}\qquad \widetilde  t^*_2:=
  \frac{\widetilde \eta_2+{\rm i}\widetilde \xi_2}{\sqrt2{\rm i}}\qquad \
   \widetilde t_3^*:=
\frac{\widetilde p+{\rm i}\widetilde q}{\sqrt2{\rm i}}\qquad\widetilde  T^*:=
  \frac{\widetilde P+{\rm i}\widetilde Q}{\sqrt2{\rm i}}  \eeqa
and we denote as
 \beqa{rps map}
\phi_{\rm rps}:\quad \textrm{\rm rps}=(\L,\widetilde\l,\widetilde\xi,\widetilde\eta,\widetilde p, \widetilde q, \widetilde P,\widetilde Q)\to (y^\ppu_{{\rm rps}}, y^\ppd_{{\rm rps}}, x^\ppu_{{\rm rps}}, x^\ppd_{{\rm rps}})
\eeqa
   respectively
       \beqa{rps complex***}
\phi^\complex_{\rm rps}:\quad \textrm{\rm rps}^\complex=(\L,\widetilde\l,\widetilde t, \widetilde t^*, \widetilde T, \widetilde T^*)\to(y^\ppu_{{\rm rps}^\complex}, y^\ppd_{{\rm rps}^\complex}, x^\ppu_{{\rm rps}^\complex}, x^\ppd_{{\rm rps}^\complex})
\eeqa
{the change of coordinates relating the real/complex rps' to the cartesian coordinates. 
The real--analytic character of $\phi_{\rm rps}$ and hence of
$\HH_{\rm rps}:=\HH\circ\phi_{\rm rps}$
on  a region of phase space having the form $\widetilde\cM_{\varepsilon_0}\times \real^2$, with $\widetilde\cM_{\varepsilon_0}:=\widetilde\cL\times\torus^2\times B^{6}_{\varepsilon_0}(0)$ has been discussed in Refs.\cite{pinzari-th09, chierchiaPi11b}. }
\\
In this paper, we are interested to look at the so--called ``retrograde problem'', whose motions live in a neighborhood of the singular manifold \beq{singularities1}\cM_\p:=\big\{\GG_1=\L_1 \text{ or } \GG_2=\L_2 \text{ or } \GG=\GG_1-\GG_2{\text{ or }\ZZ=\GG}\big\}\ .\eeq
It turns out that (see Proposition \ref{prop: analyticity domain}) the one can find regular  coordinates on a region of phase space $\cM_{\varepsilon_0}\times \real^2$, where
\beqa{M0}\cM_{\varepsilon_0}:=\cL\times \torus^2\times B^{6}_{\varepsilon_0}(0)\ .\eeqa
 including a neighborhood of $\cM_\p$ and that the formulae of such regularizing coordinates are very simply related to the ones of the prograde case, as now we briefly discuss.
}
 {We denote  as
\beqa{real and complex}
&& \textrm{\rm rps}_\p:=(\L,\l,   \eta,  \xi,  p,  q,  P, Q)=(\L_1,\L_2,\l_1,\l_2,   \eta_1,\eta_2,  \xi_1,\xi_2,  p, q,  P,  Q)\nonumber\\
&& \textrm{\rm rps}_\p^{\complex}:=(\L,\l,   t,  t^*,  T,  T^*)=(\L_1,\L_2,\l_1,\l_2,   \t_1,\t_2,  \t_3, t_1^*,  t_2^*, t_3^*, T,  T^*)\eeqa
}
 the real/complex coordinates which are defined via the formulae
\beqa{PR}
\arr{
 \L_1=\L_1\\
 \L_2=\L_2\\
t_1=\sqrt{\L_1-\GG_1}\,e^{{\rm i}(\g_1+\g+\zeta)}\\
t_2=-\ii\sqrt{\L_2-\GG_2}\,e^{{\rm i}(-\g_2+\g+\zeta)}\\
 t_3=-\ii\sqrt{\GG-\GG_1+\GG_2}\,e^{{\rm i}(\g+\zeta)}\\
T=\sqrt{\GG-\ZZ}\,e^{{\rm i}\zeta}
}\qquad 
\arr{
\l_1=\ell_1+\g_1+\g+\zeta\\
 \l_2=\ell_2+\g_2-\g-\zeta\\
t_1^*=-{\rm i}\sqrt{\L_1-\GG_1}\,e^{-{\rm i}(\g_1+\g+\zeta)}\\
 t_2^*=-\sqrt{\L_2-\GG_2}\,e^{-{\rm i}(-\g_2+\g+\zeta)}\\
t_3^*=-\sqrt{\GG-\GG_1+\GG_2}\,e^{-{\rm i}(\g+\zeta)}\\
T^*=-{\rm i}\sqrt{\GG-\ZZ}\,e^{-{\rm i}\zeta}
}
\eeqa
 and 
    \beqa{PR1}
t_1&:=&
  \frac{\eta_1-{\rm i}\xi_1}{\sqrt2}\qquad t_2:=
  \frac{\ii\eta_2-\xi_2}{\sqrt2}\qquad\ \ \ t_3:=
  \frac{\ii p-q}{\sqrt2}\qquad\ \ \ T:=
  \frac{P-{\rm i}Q}{\sqrt2}\nonumber\\\nonumber\\
  t^*_1&:=&
  \frac{\eta_1+{\rm i}\xi_1}{\sqrt2{\rm i}}\qquad  t^*_2:=
  \frac{\ii\eta_2+\xi_2}{\sqrt2{\rm i}}\qquad \
   t_3^*:=
\frac{\ii p+q}{\sqrt2{\rm i}}\qquad T^*:=
  \frac{P+{\rm i}Q}{\sqrt2{\rm i}}.  \eeqa
  Observe that $\big\{(t_1, t_2, t_3, T,t_1^*, t_2^*, t_3^*, T^*)=(0,0,0,0, 0,0,0,0)\big\}$ is the equation of the manifolds $\cM_\p$'s.
{We let
       \beqa{rps+ map}
&&\phi_{{\rm rps}_\p}:\quad {{\rm rps}_\p}=(\L,\l,\eta, \xi, p,  q,  P, Q)\to (y_{{\rm rps}_\p}^\ppu, y_{{\rm rps}_\p}^\ppd, x_{{\rm rps}_\p}^\ppu, x_{{\rm rps}_\p}^\ppd)\nonumber\\
&&\phi_{{\rm rps}^\complex_\p}:\quad {{\rm rps}^\complex_\p}=(\L,\l, t,t^*, T,T^*)\to (y_{{\rm rps}^\complex_\p}^\ppu, y_{{\rm rps}^\complex_\p}^\ppd, x^\ppu_{{\rm rps}^\complex_\p}, x^\ppd_{{\rm rps}^\complex_\p})
\eeqa}
Let $\phi_2^-$ be  the involution 
\beqa{involution}\phi_2^-\big(\L_1,\L_2,\l_1, \l_2,t,t^*, T, T^*\big):=\big(\L_1,-\L_2,\l_1, -\l_2,t,t^*, T, T^*\big)\ .\eeqa Then
\begin{proposition}\label{signs1} $\phi_{{\rm rps}^\complex_\p}=\phi_{{\rm rps}^\complex}\circ\phi_2^-$. \end{proposition}
Proposition \ref{signs1} (the proof of which can be found in Appendix~\ref{lemma sui segni}) is useful in order to derive the explicit expression of $\phi_{\rm rps_\p}$ from the one of $\phi_{\rm rps_\p}$, as done in Appendix \ref{formulae of reversed map}. 
The latter is needed to determine
a  analyticity domain for the Hamiltonian 
\beqa{rps Ham real}
\HH_{{\rm rps}_\p}&=&\HH\circ\phi
 _{{\rm rps}_\p}
 =-\frac{\mm_1^3\MM_1^2}{2\L_1^2}-\frac{\mm_2^3\MM_2^2}{2\L_2^2}+\m\, \Big(
 \frac{y_{{\rm rps}_\p}^\ppu\cdot y_{{\rm rps}_\p}^\ppd}{m_0}-\frac{m_1 m_2}{|x_{{\rm rps}_\p}^\ppu-x_{{\rm rps}_\p}^\ppd|}
 \Big)\nonumber\\
 &=:&\hh_{\rm k}(\L)+\m f_{\rm rps_\p}(\L,\l, \eta,\xi,p,q)\eeqa
 which we shall provide in the next section.

\subsubsection{\label{analyticity domain}A domain of analyticity for $\HH_{\rm rps_\p}$}
{
Let $\chi$,  $\a_-$, $\a_+$,   be pure numbers verifying
\beqa{domain rps+}\chi>1\ ,\quad  0<\a_-<\a_+<1\ .\eeqa
Let the masses $m_0$, $ m_1$, $ m_2$ and the number $\m$ in~\equ{helio}--\equ{reduced masses} be chosen so that
\beqa{domain rps+1} m_2<\frac{\sqrt{\a_-}}{2\chi}m_1\ ,\quad 0<\m<\frac{3m_0}{m_1-4 m_2}\eeqa
and let,  for fixed $\L_-$, $\L_+$, with \beqa{domain rps+2}0<\L_-<\L_+\ ,\eeqa
  \beq{L0}\cL:=\Big\{\L=(\L_1,\L_2): \ \L_-< \L_2< \L_+\ ,\quad k_-\L_2< \L_1< k_+\L_2\Big\}\eeq
where
 \beq{kpm}k_\pm:=\frac{m_1}{m_2}\sqrt{\frac{m_0+\m m_2}{m_0+\m m_1}\a_\pm}\ .\eeq
 We the take, for the coordinates ${\rm rps}_\p=(\L,\l,\eta,\xi,p,q,P,Q)$ in~\equ{rps+ map} the domain $\cM_{\varepsilon_0}\times\Big\{(0,0)\Big\}$, where $\cM_{\varepsilon_0}$ is as in \equ{M0}.
 \begin{remark}\rm
 \item[$\star$] The cyclic pair $(P,Q)$ has been fixed, as in~Refs.\cite{pinzari-th09, chierchiaPi11b},  to $(0,0)$. This corresponds to fix the direction of the third axis of the  prefixed reference frame parallel to the total angular momentum direction and has not influence on the dynamics.
 \item[$\star$] The bound in~\equ{domain rps+1} for $m_2$ is not purely technical, but necessary in order to realize the {\sc orc} configuration, with arbitrarily small eccentricities and inclinations. Indeed, when the motions are planar, circular and {\sc orc}, one has $\L_1=|\CC^\ppu|>|\CC^\ppd|=\L_2$, which implies, as discussed above,~\equ{domain rps+1}.
\item[$\star$]
Recalling that (by~\equ{reduced masses} and~\equ{aL})  the semi--axes ratio $\a:=\frac{a_1}{a_2}$ is related to the ratio $\frac{\L_1}{\L_2}$ via
 \beqa{alpha ratio}\frac{\L_1}{\L_2}&=&\frac{{\rm  m}_1}{{\rm  m}_2}\sqrt{\frac{{\rm  M}_1}{{\rm  M}_2}\a}=\frac{m_1}{m_2}\sqrt{\frac{m_0+\m m_2}{m_0+\m m_1}\a}\eeqa
and that, by definition of $ \cL$, $k_-<\frac{\L_1}{\L_2}<k_+$ and using the definitions of $k_\pm$ in~\equ{kpm}, we have
\beqano\a_-<\a<\a_+\ .\eeqano
Below, we shall choice $\a_+<\a_s$, where
$\a_s$, $\cL_s$ are as in~\equ{Ls}, so that we shall be able to use the results of the previous sections.
\end{remark}
 \begin{proposition}\label{prop: analyticity domain}
One can find $\varepsilon_0>0$, depending only on $\L_-$, $\chi$, $\a_+$ such that
the function $\HH_{\rm rps_\p}$ in \equ{rps Ham real}
 is real--analytic for $(\L,\l,\eta,\xi,p,q)\in\cM_{\varepsilon_0}$, where $\cM_0$ is as in \equ{M0}, for all choices of $m_0$, $m_1$, $m_2$, $\a_-$, $\L_+$ and $\m$ satisfying~\equ{domain rps+},~\equ{domain rps+1} and~\equ{domain rps+2}.
 \end{proposition}
We have to ensure that, on the domain $\cM_{\varepsilon_0}$, the map $\phi_{{\rm rps}_\p}$ discussed in Section \ref{formulae of reversed map} is real--analytic and collisions are excluded.
Observe  that, as a consequence of~\equ{domain rps+1}, $\frac{m_1}{m_2}\sqrt{\a_-}>2\chi$, $\sqrt{\frac{m_0+\m m_2}{m_0+\m m_1}}>\frac{1}{2}$, whence
$\chi<k_-<k_+$.
This inequality implies, on $\cL$,
$
|\L_1-\L_2|\ge \L_-(\chi-1)
$
and therefore $|\L_1-\L_2|\ge\frac{1}{2} \L_-(\chi-1)$ on a complex neighborhood of $\cL$ depending on $\L_-$, $\chi$.
Therefore, the former common denominator of ${\rm c}_1^*$, ${\rm c}_2$, ${\rm c}^*_2$ in~\equ{c} well separated from zero as soon as $\varepsilon<\frac{1}{4}\sqrt{\L_-(\chi-1)}$. The other denominators are treated in a similar way. The fact that the 
 planar Kepler maps~\equ{planar poinc} are analytic  and that collisions are excluded on a complex neighborhood of  $\cM_\varepsilon$ depending only on $\L_-$, $\a_+$ and $\chi$
are classical facts and will be not discussed.
}


\subsubsection{D'Alembert rules}
The Hamiltonian $\HH$ in~\equ{helio} is invariant by rotations of the reference frame  and reflections of the coordinates axes. 
These symmetries reflect on the Hamiltonian  $\HH_{\rm rps_\p}$ in \equ{rps Ham real} and  they will be named here ``D'Alembert rules'',  by similarity with analogue relations holding 
when $\HH$ is written in Poincar\'e coordinates, where they are given such name.
A full analytical discussion may be found in Appendix \ref{DAlembert rules}. In this section, we summarize their effect on the averaged perturbing function. More precisely, we focus on the system written in complex coordinates
\beqa{RPS Ham}
\HH_{{\rm rps}^\complex_\p}(\L,\l,t,t^*):=\HH\circ \phi_{{\rm rps}^\complex_\p}(\L,\l,t,t^*)&=&\hh_{\rm k}(\L_1,\L_2)+\m f_{{\rm rps}^\complex_\p}(\L,\l,t,t^*)
\eeqa
and we let
$${f^{\rm av}_{{\rm rps}^\complex_\p}}(\L,t, t^*):=\frac{1}{(2\p)^2}\int_{[0,2\p]^2}f_{{\rm rps}^\complex_\p}(\L,\l,t,t^*) d\l\ .$$
Let $c_{a,a^*}$ be the coefficients of the expansion of $f_{\rm rps^\complex_\p}^{\rm av}$ up to a prefixed order $m$
\beqa{retr exp}{ f^{\rm av}_{\rm rps^\complex_\p}}=\sum_{|(a,a^*)|\le m} c_{a,a^*}(\L) t^a{ {t^*}}^{a^*}+{\rm o}_m(t,t^*;\L)
\eeqa
(which are well defined if the corresponding real coordinates vary in the domain $\cM_{\varepsilon_0}$ as in Proposition \ref{prop: analyticity domain}) then
\begin{itemize}
\item[$({\rm s}_1)$] by~\equ{orc DAlembert},
 $ c_{a,a^*}(\L)\ne 0\ \Longrightarrow\ \sum_{i}a_i=\sum_i a_i^*$;
\item[$({\rm s}_2)$] by the latter identity in~\equ{rot and reflect*} and previous item, ${ f^{\rm av}_{\rm rps_\p}}$ is even in $(t_3,t_3^*)$.  Therefore, by the previous item, ${ f^{\rm av}_{\rm rps_\p}}$ is even in $(t_3,t_3^*)$ and $(t_1,t_2,t_1^*,t_2^*)$ separately;
\item[$({\rm s}_3)$]  by the former identity in~\equ{rot and reflect*}  and the parity in $(t_1,t_2,t_1^*,t_2^*)$, the coefficients $ c_{a,a^*}$ in~\equ{retr exp}  do not change letting   $$(a_1,a_2,a_1^*,a_2^*)\to (a_1^*,a_2^*, a_1,a_2)$$
  and leaving $(a_3, a_3^*)$ unvaried;
  \item[$({\rm s}_4)$] by the reality conditions ${ f^{\rm av}_{\rm rps_\p}}=({ f^{\rm av}_{\rm rps_\p}})^{\rm cc}$, $$c_{a,a^*}=(-1)^{\frac{|a|+|a_*|}{2}} (-1)^{a_2+a_2^*+a_3+a_3^*} c_{a^*, a}$$
as it follows using
the following relations (implied by \equ{PR1})
\beqa{sym retrogr}\Big( t_1,  t^*_1\Big)^{\rm cc}=\ii \Big( t_1^*,  t_1\Big)\ ,\quad \Big( t_j,  t^*_j\Big)^{\rm cc}=-\ii \Big( t_j^*,  t_j\Big)
\ ,\quad j=2,\ 3
\eeqa
into \equ{retr exp}.
   \end{itemize} 
\subsubsection{\label{PRBNF}Birkhoff theory in the  partially reduced retrograde problem}
{A remarkable consequence of $({\rm s}_1)$--$({\rm s}_4)$ is that (as already mentioned in the introduction), differently from what happens in the prograde problem, in the case of the retrograde one, it is {\it not} possible to guarantee a priori that the eigenvalues of the matrix $\s$ appearing in the second--order  expansion
\beqa{quadratic retrograde expr}
f^{\rm av}_{{\rm rps}^\complex_\p}=C_0(\L_1,\L_2)+\ii t_h\cdot \s(\L) t^*+\ii\varsigma(\L)t_3 t_3^*+{\rm O}_4(t,t^*;\L)\eeqa
where $t_h:=(t_1,t_2)$, $t^*_h:=(t^*_1,t^*_2)$, are real. 
Indeed, the ${\rm s}_1$--${\rm s}_4$  imply that the entries $\s_{ij}$ of $\s$ and $\varsigma$  in \equ{quadratic retrograde expr} verify
\beqa{reals}\s_{11} \ \s_{22}\in \real\ ,\quad \s_{12}=\s_{21}\in \ii\real\ ,\quad \varsigma\in \real\eeqa
(in the case of the prograde problem, studied in Refs\cite{pinzari-th09, chierchiaPi11b}, one has, correspondingly, that  $\s$ is symmetric and real, because of the different form of the reality condition ${\rm s}_4$).
}
The following result is therefore not trivial:

\begin{proposition} \label{non res orc} {\bf(Herman resonance, Birkhoff theory and symmetries in the partially reduced retrograde problem)} \\
Let $\cL$ be as in \equ{L0}, $\varepsilon_0$ as in Proposition \ref{prop: analyticity domain}. 
\item[{\rm (i)}]
For $\L\in \cL$,
the matrix
$\s(\L)$ in~\equ{quadratic retrograde expr} has two distinct eigenvalues $\s_1(\L)$, $\s_2(\L)$ which, together with $\varsigma$, are
 real--analytic  for all $\L\in\cL$. 
Moreover,  
for all $\L\in \cL$,
one can find a 
a real--analytic vector matrix $\O:\ \L\in \cL\to \O(\L)\in \real^3$, a
symplectic, real--analytic  $4\times 4$ function matrix $V:\ \L\in \cL\to V(\L)$   and a   real--analytic, canonical transformation
\beqa{real phi}
\check\phi:\quad\check\cM_{\check\varepsilon}:= \cL\times \torus^2\times B^{6}_{\check\varepsilon}(0)&\to&\cM_{\varepsilon_0}\nonumber\\
(\L,\check\l, (\check \eta, \check\xi), \check p,  \check q)&\to&(\L,\check\l+\check\f(\L,\check \eta, \check \xi), V(\L) (\check \eta, \check\xi),  \check p,  \check q)
\eeqa
such that, if \beqa{diagonalized orc complex}
\check\HH:=\HH\circ\check\phi=\hh_{\rm k}(\L)+\m\check f(\L,\check\l,\check \eta, \check \xi, \check p, \check q)
\eeqa
 then   \beqa{diagonalized complex}
\check f^{\rm av}=C_0(\L)+
\Omega(\L)\cdot\check\t
+{\rm O}_4(\check\eta,\check\xi,\check p, \check q;\L)\quad \eeqa
with $$\check\t=\left(\frac{\check\eta_1^2+\check\xi_1^2}{2},\ \frac{\check\eta_2^2+\check\xi_2^2}{2}\ ,\ \frac{\check p^2+\check q^2}{2}\right)$$
Moreover, for any prefixed $s\in \natural$, one can find 
$\a_s$ such that, if
\beqa{Ls}\cL_s:=\big\{\L\in{\cal L}:\quad 0<\a<\a_s\big\}\eeqa
\beqa{non res prograde}\forall\ \L\in \cL_s\ ,\ k\in \integer^3\ ,\quad0<|k|\le 2s\ ,\quad  k\ne m(1,1,1)\ ,\quad m\in \integer:\quad 
\O(\L)\cdot k\ne 0
\ .\eeqa
\item[{\rm (ii)}]
It is possible to find 
a real--analytic $3\times 3$ matrix function ${\rm T}:\ \L\in \cL_s\to {\rm T}(\L)$ and real--analytic functions $(\L,\ovl\t)\in \cL_s\times \complex^3\to \cP_j(\L,\ovl\t)$ with $ \cP_j(\L,\ovl\t)$ a polynomial of degree $j=3$, $\cdots$, $s$ in $\ovl\t$
and a real--analytic and canonical transformation
\beqano
\ovl\phi:\quad\ovl\cM_{\ovl\varepsilon}:= \cL_{s}\times \torus^2\times B^{6}_{\ovl\varepsilon}(0)&\to& \check\cM_{\check\varepsilon}=\cL_s\times \torus^2\times B^{6}_{\check\varepsilon}(0)\nonumber\\
(\L,\ovl\l, \ovl \eta, \ovl \xi, \ovl p, \ovl q)&\to&(\L,\check\l, \check \eta, \check \xi, \check p, \check q)=(\L,\ovl\l+\ovl\f(\L,\ovl \eta, \ovl \xi, \ovl p, \ovl q), \ovl B(\ovl \eta, \ovl \xi, \ovl p, \ovl q;\L))
\eeqano
such that, if
\beqa{normalized orc}
\ovl\HH:=\check\HH\circ\ovl\phi=\hh_{\rm k}(\L)+\m\ovl f(\L,\ovl\l,\ovl \eta, \ovl \xi, \ovl p, \ovl q)\quad (\L,\ovl\l,\ovl \eta, \ovl \xi, \ovl p, \ovl q)\in \cL_{2s}\times \torus^2\times B^{6}_{\ovl\varepsilon}(0)
\eeqa
then
\beqa{BNF retrograde}
\ovl f^{\rm av}
=C_0(\L)+\O\cdot\ovl\t+\frac{1}{2}\ovl\t\cdot{\rm T}(\L)\ovl\t+\sum_{j=3}^s\ovl\cP(\ovl\t;\L)+{\rm O}_{2s+2}(\ovl\eta,\ovl\xi, \ovl p, \ovl q;\L)
\eeqa
where $$\ovl\t=\left(\frac{\ovl\eta_1^2+\ovl\xi_1^2}{2},\ \frac{\ovl\eta_2^2+\ovl\xi_2^2}{2}\ ,\ \frac{\ovl p^2+\ovl q^2}{2}\right)\ ;$$
\item[{\rm (iii)}] If $\widetilde\O(\L_1,\L_2)$ and $\widetilde{\rm T}(\L_1,\L_2)$ denote the first and second order Birkhoff invariants for the prograde problem with $n=2$ (denoted in Ref.\cite[Proposition 7.1 and Equation (8.1)]{chierchiaPi11b}, respectively, as $(\s, \ovl\varsigma)$, $\t$) then 
$\O$, ${\rm T}$ are related to $\widetilde\O$, $\widetilde{\rm T}$ in via
 \beqa{TT}
\O_i(\L_1,\L_2)=s_i\widetilde\O_i(\L_1,-\L_2)\ ,\quad {\rm T}_{ij}(\L_1,\L_2)=s_is_j\widetilde{\rm T}_{ij}(\L_1,-\L_2)\ 1\le i\le j\le 3
\eeqa
with $s_1=-s_2=-s_3=1$;
\item[{\rm (iv)}]
$\check\phi$ and $\breve\phi$ commute with $\cR_g$, $\cR_3^-$ and $\cR_{1\leftrightarrow2}$ in~\equ{Rg} and~\equ{orc DAlembert1}.
\end{proposition}
\begin{remark}\rm
{The transformation $\check\phi$ preserves the function 
\beqa{G}\GG
=\L_1-\L_2-\frac{\eta_1^2+\xi_1^2}{2}+\frac{\eta_2^2+\xi_2^2}{2}+\frac{p^2+q^2}{2}
\eeqa
  Differently from the prograde case, $V(\L)\notin{\rm SO(4)}$.
 }

\end{remark}
 \vskip.1in
 {\bf Proof.} 
(i)  The transformation $\tilde\phi$ in~Ref.\cite[Equations (7.20) and (7.22)]{chierchiaPi11b} 
  will be here denoted at
\beqa{hat U}\widehat\phi:\ (\L,\widehat\l, \widehat z)\to (\L, \widetilde\l, \widetilde z)=(\L,\widehat\l+\widehat\psi(\L, \widehat z), \widehat U(\L)\widehat z)\eeqa
  with $\widehat z=(\widehat\eta, \widehat\xi, \widehat p, \widehat q)$, $\widetilde z=(\widetilde\eta, \widetilde\xi, \widetilde p, \widetilde q)$.
$\widehat\phi$
 projects as a transformation from 
 \beqano
\widehat\phi^\complex:\quad \cL\times \torus^2\times B^{6}_{\widehat\varepsilon}(0)&\to& \cL\times \torus^2\times B^{6}_{\widehat\varepsilon}(0)\nonumber\\
(\L,\widehat\l, (\widehat t, \widehat t^*))&\to&(\L,\widetilde\l, \widetilde t, \widetilde  t^*)=(\L,\widehat\l+\widehat\psi(\L,\widehat t, \widehat t^*), (\widehat U\widehat t, \widehat U\widehat t^*))
\eeqano
where$(\widetilde t, \widetilde t^*)$,   $(\widehat t, \widehat t^*)$ are related to $\widetilde z=(\widetilde\eta,\widetilde\xi, \widetilde p, \widetilde q)$, $\widehat z=(\widehat\eta,\widehat\xi, \widehat p, \widehat q)$ at the same way as in~\equ{PROld+}, respectively.
We then take $\check\phi^\complex:=\phi_2^-\circ\widehat\phi^\complex\circ\phi_2^-$, where $\phi_2^-$ in as in~\equ{involution}. The transformation $\check\phi^\complex$ is 
\beqa{check psi}
\check\phi^\complex:\quad \cL\times \torus^2\times B^{6}_{\check\varepsilon}(0)&\to& \cL\times \torus^2\times B^{6}_{\check\varepsilon}(0)\nonumber\\
(\L,\check\l, \check t, \check t^*)&\to&(\L,\l,  t,  t^*)=(\L,\check\l+\check\psi(\L,\check t, \check t^*), (U\check t,  U\check t^*))
\eeqa
 with 
$\check\psi_1=\widehat\psi_1(\L_1,-\L_2,\check t, \check t^*)$,
$\check\psi_2=-\widehat\psi_2(\L_1,-\L_2,\check t, \check t^*)$, $U(\L_1,\L_2)=\widehat U(\L_1,-\L_2)$. Then
$U$ is symmetric and verifies \beqa{eqU}U^{-1}(\L) \s(\L) U(\L)=\diag(\s_1(\L),\s_2(\L))=:D(\L)\eeqa 
where $\s_1,\s_2$ are the eigenvalues of $\s$, with $\s_i(\L_1,\L_2)=\widetilde\s_i(\L_1,-\L_2)$. Observe that, while $U$ needs not to be real, since $\widehat U\in {\rm SO}(2)$,  $U$ verifies
\beq{U is unique}U^{\rm t}=U^{-1}\ ,\quad \det U=1\ .\eeq
Note that $U$ is uniquely determined by conditions~\equ{eqU} and~\equ{U is unique}.
Let $(\L,\check\l, \check\eta, \check\xi, \check p, \check q)$ be related to $(\L,\check\l, \check t, \check t^*)$ via~\equ{PR1}, and let $\check\phi$ be corresponding projection $\check\phi^\complex$ as a transformation from $(\L,\check\l, \check\eta, \check\xi, \check p, \check q)$ to $(\L,\l, \eta, \xi,  p, q)$ of the form in~\equ{real phi}, with
\beq{phi V}\check\f(\L,\check\eta, \check\xi)=\check\psi(\L, L(\check\eta, \check\xi))\quad V(\L)=L^{-1}\diag(U(\L),U(\L))L\ ,\eeq
where
\beqa{L}L=
\left(
\begin{array}{ccccc}
\frac{1}{\sqrt2}&0&-\frac{\ii}{\sqrt2}&0\\
0&\frac{\ii}{\sqrt2}&0&-\frac{1}{\sqrt2}\\
-\frac{\ii}{\sqrt2}&0&\frac{1}{\sqrt2}&0\\
0&\frac{1}{\sqrt2}&0&-\frac{\ii}{\sqrt2}
\end{array}
\right)\ ,\qquad L^{-1}=
\left(
\begin{array}{ccccc}
\frac{1}{\sqrt2}&0&\frac{\ii}{\sqrt2}&0\\
0&-\frac{\ii}{\sqrt2}&0&\frac{1}{\sqrt2}\\
\frac{\ii}{\sqrt2}&0&\frac{1}{\sqrt2}&0\\
0&-\frac{1}{\sqrt2}&0&\frac{\ii}{\sqrt2}
\end{array}
\right)\ .
\eeqa

By construction, 
using
also
$$\frac{\check\eta_1^2+\check\xi_1^2}{2}=\ii \check t_1\check t_1^*\ ,\quad -\frac{\check\eta_2^2+\check\xi_2^2}{2}=\ii \check t_2\check t_2^*\ ,\quad -\frac{\check p^2+\check q^2}{2}=\ii \check t_3\check t_3^*$$
$\check\phi$ verifies~\equ{diagonalized orc complex} and~\equ{diagonalized complex}, with  $\sigma_1(\L_1,\L_2)=\widetilde\sigma_1(\L_1,-\L_2)$, $\sigma_2(\L_1,\L_2)=\widetilde\sigma(\L_1,-\L_2)$, $\varsigma(\L_1,\L_2)=\widetilde\varsigma(\L_1,-\L_2)$. We have to check that $\s_1$, $\s_2$, $\varsigma$, $\check\f$ and $V$ are real when their arguments are so.  It follows from~Ref.\cite[Equations (B.1) and (B.2)]{chierchiaPi11b} and Proposition \ref{signs1} that $f^{\rm av}_{\rm rps^\complex_\p}$ has the expression in~\equ{quadratic retrograde expr}, with $C(\L_1,\L_2)=\widetilde C(\L_1,-\L_2)$ and
\beqa{quadratic retrograde}  \s(\L_1,\L_2)=\widetilde\s(\L_1,-\L_2)\ ,\quad \varsigma(\L_1,\L_2)=\widetilde\varsigma(\L_1,-\L_2)\ .\eeqa
Explicitly, using Ref.\cite[Equation (B.1)]{chierchiaPi11b}, we find (according to \equ{reals})
\beqano\s(\L_1,\L_2)=\left(
  \begin{array}{ccc}
\frac{\rm  s}{\L_1}&-{\rm i}\frac{\widetilde{\rm  s}}{\sqrt{\L_1\L_2}}\\
-{\rm i}\frac{\widetilde{\rm  s}}{\sqrt{\L_1\L_2}}&-\frac{\rm  s}{\L_2}
  \end{array}
\right)\qquad\varsigma(\L)=-\Big(\frac{1}{\L_1}-\frac{1}{\L_2}\Big){\rm  s}\ ,\eeqano
with
 \beqano
 {\rm  s}:=- m_1m_2\frac{\a}{2a_2}{b^{(1)}_{3/2}}(\a)\qquad \widetilde{\rm  s}:=m_1m_2\frac{\a}{2a_2}{b^{(2)}_{3/2}(\a)}
\eeqano
with $\a$, $a_2$ as in~\equ{aL},~\equ{alpha ratio},
and, as usual,
$b^{(j)}_s(\a)$'s being the Laplace coefficients, defined via the Fourier expansion $$\frac{1}{\big(1-2\a\cos\theta+\a^2\big)^s}=\sum_{k\in \integer}b^{(k)}_s(\a)e^{{\rm i} k\theta}\ \qquad {\rm i}:=\sqrt{(-1)}\ .$$
 The  eigenvalues of $\s$ can be explicitly computed and they turn to be real. Indeed,
\beq{elleq1}\s_{1}, \s_2=\frac{\tr\s}{2}\pm\frac{1}{2}\sqrt{(\tr\s)^2-4\det\s}\ .\eeq
{Since $\tr\s=\Big(\frac{1}{\L_1}-\frac{1}{\L_2}\Big){\rm  s}$ is real, we} have  to check that the discriminant 
\beqano
\D:=(\tr\s)^2-4\det\s=(\frac{1}{\L_1}-\frac{1}{\L_2})^2{\rm s}^2+\frac{4}{\L_1\L_2}\big({\rm  s}^2-\widetilde{\rm  s}^2\big)
\eeqano
is positive.
Recalling that the Laplace coefficients  verify
$$b^{(j)}_s(\b)> b^{(j+1)}_s(\b)\quad \textrm{for all}\quad s>0,\quad j\in \integer,\quad 0<|\b|<1,$$
(see~Ref.\cite{fejoz04} for a proof), one has
\beq{elleq2}
{\rm  s}^2-\widetilde{\rm  s}^2=(m_1m_2\frac{\a}{a_2})^2\big((b^{(1)}_{3/2}(\a))^2-(b^{(2)}_{3/2}(\a))^2\big)> 0.\eeq
and we have the assertion. 

To prove the reality of $V(\L)$, we first  check that the matrix $U(\L)$ 
has the form\
\beqa{hyperbolic matrix}U(\L)=\left(
\begin{array}{ccc}
\cosh x(\L)&-\ii\sinh x(\L)\\
\ii\sinh x(\L)&\cosh x(\L)
\end{array}
\right)\eeqa
with some $x=x(\L)\in \real$.
Let, for a given $2\times 2$ matrix $A$, $\cC(A)$ be defined via $\cC(A)_{ii}=A_{ii}^{\rm cc}$; $\cC(A)_{ij}=-A_{ii}^{\rm cc}$ for $i\ne j$. Since $\cC(\s)=\s$, $\cC(D)=D$ and $\cC(AB)=\cC(A)\cC(B)$, applying $\cC$ to the equality~\equ{eqU}, we find $\cC(U)=U$, which is equivalent to

 \beq{U}U_{ii}(\L) \in \real\ ,\quad U_{ij}(\L)\in \ii\real\quad \forall \ i\ne j=1,\ 2\ .\eeq
Taking in count~\equ{U is unique},~\equ{U},  we find that $U$ has necessarily the form~\equ{hyperbolic matrix}.
Then, by direct computation, we find that $V(\L)$ in~\equ{phi V} is given by 
\beqa{V}V(\L)=\diag(W(x(\L)),W(-x(\L)))\quad {\rm with}\quad W(x)=\left(
\begin{array}{ccc}
\cosh x& \sinh x\\
\sinh x&\cosh x
\end{array}
\right)\ ,\eeqa
We now check that $\f(\L,\check\eta,\check\xi)$
is real.
We use the generating function of $\check\phi^\complex$, which is
$$S(\check\L, \l, \check t, t^*)=\check\L \l+t_h^*\cdot U\check t_h+\check t_3 t_3^*\ .$$
Then the function $\check\psi$ in~\equ{check psi} has the analytical form
$$\psi(\L,\check t_h, \check t_h^*)=\check t_h^*\cdot A(\L)\check t_h\quad {\rm where}\quad A(\L):=U(\L)^{-1}\partial_\L U(\L) \ .$$ 
with
 $A(\L)$  skew--symmetric, because of~\equ{U is unique}. Using~\equ{U} we find
$$A(\L)=\partial_\L x(\L)\left(
\begin{array}{ccc}
0&-\ii\\
\ii&0
\end{array}
\right)\ ,$$
which gives
\beqa{psi}\check\psi(\L,\check t_h, \check t_h^*)=\ii\partial_\L x(\L)\big(\check t_1\check t_2^*-\check t_1^*\check t_2\big)\ .\eeqa
Then, by the first equation in~\equ{phi V} and~\equ{L}, 
$$\check\f(\L,\check\eta, \check\xi)=\partial_\L x(\L)\big( \check\eta_1 \check\xi_2+ \check\xi_1 \check\eta_2\big)$$
is real on real arguments.
The proof of \equ{non res prograde} goes as in Ref.\cite[Proof of Proposition 7.2]{chierchiaPi11b} upon replacing $\s_1$, $\s_2$ $\varsigma$ with
$$\s_1=-\frac{3}{4\L_1}\frac{a^2_1}{a_3^2}+{\rm O}(\frac{a_1^3}{a_2^4\L_1})\qquad \s_2=+\frac{3}{4\L_2}\frac{a_1^2}{a_2^3}+{\rm O}(\frac{a_1^2}{a_2^3\L_1})\ ,\quad \varsigma=\frac{3}{4}\frac{a^2_1}{a_3^2}+{\rm O}(\frac{a_1^3}{a_2^4\L_1})\left(\frac{1}{\L_1}-\frac{1}{\L_2}\right)\ .$$
(ii)
We denote as
\beqano
\breve\phi:\quad \widetilde\cL\times \torus^2\times B^{6}_{\breve\varepsilon}(0)&\to& \widetilde\cL\times \torus^2\times B^{6}_{\breve\varepsilon}(0)\nonumber\\
(\L,\breve\l, \breve z)&\to&(\L,\widehat\l, \widehat z)=(\L,\breve\l+\breve\psi(\L,\breve z), \breve B'(\L,\breve z))
\eeqano
with $\breve z=(\breve\eta, \breve\xi, \breve p, \breve q)$m $\widehat z=(\widehat\eta, \widehat\xi, \widehat p, \widehat q)$
the
 transformation  in Ref.\cite[Equation (7.28)]{chierchiaPi11b}.
 As in the proof of Proposition 	\ref{non res orc}, we denote as $(\L,\widehat\l,\widehat t, \widehat t^*)$, $(\L,\breve\l,\breve t, \breve t^*)$ the coordinates related to $(\L,\widehat\l,\widehat z)$, $(\L,\breve\l,\breve z)$ via the relations in~\equ{PROld+}, respectively, and we let
 \beqano
\breve\phi^\complex:\quad \widetilde\cL\times \torus^2\times B^{6}_{\breve\varepsilon}(0)&\to& \widetilde\cL\times \torus^2\times B^{6}_{\breve\varepsilon}(0)\nonumber\\
(\L,\breve\l, \breve \eta, \breve t,\breve t^*)&\to&(\L,\widehat\l, \widehat t,  \widehat t^*)=(\L,\breve\l+\breve\psi(\L,\breve \eta, t,\breve t^*), \breve B'(\breve t,\breve t^*;\L))
\eeqano
the transformation induced by $\breve\phi$ on such  coordinates. Then, by definition,
\beqano
\breve\HH^\complex:=\HH\circ\breve\phi^\complex=\hh_{\rm k}(\L)+\m\breve f^\complex(\L,\breve\l,\breve t, \breve t^*)
\eeqano
with
\beqano
{\breve f^\complex}^{\rm av}=C_0(\L)+\widetilde\O\cdot\breve\t+\frac{1}{2}\breve\t\cdot\breve{\rm T}(\L)\breve\t+\sum_{j=3}^s\breve\cP(\breve\t;\L)+\breve\cR^\complex_{2s+2}(\L,\breve t, \breve t^*)
\eeqano
where
$$\breve\t_1=\frac{\breve\eta_1^2+\breve\xi_1^2}{2}=\ii \breve t_1\breve t_1^*\ ,\quad \breve\t_2=\frac{\breve\eta_2^2+\breve\xi_2^2}{2}=\ii \breve t_2\breve t_2^*\ ,\quad \breve\t_3=\frac{\breve p^2+\breve q_1^2}{2}=\ii \breve t_3\breve t_3^*$$
If $\phi_2^-$ is as in~\equ{involution}, we let
$\ovl\phi^\complex:=\phi_2^-\circ \breve\phi^\complex\circ \phi^\complex$ and
$\ovl\HH^\complex:=\check\HH\circ\ovl\phi^\complex$. 

Then $\ovl\phi^\complex$ and $\ovl\HH$ have the form
\beqano
\ovl\phi^\complex:\quad &&\cL_{s}\times \torus^2\times B^{6}_{\ovl\varepsilon}(0)\to \cL_s\times \torus^2\times B^{6}_{\ovl\varepsilon}(0)\nonumber\\
&&(\L,\ovl\l, \ovl \eta, \ovl t,\ovl t^*)\to(\L,\check\l, \check t,  \check t^*)=(\L,\ovl\l+\ovl\psi(\L,\ovl \eta, t,\ovl t^*), \ovl B'(\ovl t,\ovl t^*;\L))\nonumber\\
&&
\ovl\HH^\complex=\HH\circ\ovl\phi^\complex=\hh_{\rm k}(\L)+\m\ovl f^\complex(\L,\ovl\l,\ovl t, \ovl t^*)
\eeqano
with
\beqano
{\ovl f^\complex}^{\rm av}=\ovl C_0(\L)+\ovl\O\cdot\ovl\t+\frac{1}{2}\ovl\t\cdot\ovl{\rm T}(\L)\ovl\t+\sum_{j=3}^s\ovl\cP(\ovl\t;\L)+\ovl\cR^\complex_{2s+2}(\L,\ovl t, \ovl t^*)
\eeqano
where $\ovl C_0(\L_1,\L_2)=C_0(\L_1,-\L_2)$, $\ovl\O(\L_1,\L_2)=\widetilde\O(\L_1,-\L_2)=(\s_1,\s_2,\varsigma)$,
$$\ovl\t_1=\ii \ovl t_1\ovl t_1^*\ ,\quad \ovl\t_2=\ii \ovl t_2\ovl t_2^*\ ,\quad \ovl\t_3=-\ii \ovl t_3\ovl t_3^*$$
and $\ovl{\rm T}(\L_1,\L_2)=\breve{\rm T}(\L_1,-\L_2)$.
Now, if  $(\L,\ovl\l,\ovl\eta,\ovl\xi,\ovl p, \ovl q)$ are related to $(\L,\ovl\l,\ovl t, \ovl t^*)$   via the relations in~\equ{PR1},
 $\ovl\phi^\complex$, induces a transformation from $(\L,\ovl\l,\ovl\eta,\ovl\xi,\ovl p, \ovl q)$ to $(\L,\check\l,\check\eta,\check\xi,\check p, \check q)$ having the form
\beqano
\ovl\phi:\quad &&\cL_{s}\times \torus^2\times B^{6}_{\ovl\varepsilon}(0)\to \cL_s\times \torus^2\times B^{6}_{\ovl\varepsilon}(0)\nonumber\\
&&(\L,\ovl\l,\ovl\eta,\ovl\xi,\ovl p, \ovl q)\to(\L,\check\l,\check\eta,\check\xi,\check p, \check q)=(\L,\ovl\l+\ovl\f(\L,\ovl\eta,\ovl\xi,\ovl p, \ovl q), \ovl B(\ovl\eta,\ovl\xi,\ovl p, \ovl q;\L))\nonumber\\
&&
\ovl\HH=\HH\circ\ovl\phi=\hh_{\rm k}(\L)+\m\ovl f(\L,\ovl\eta,\ovl\xi,\ovl p, \ovl q)
\eeqano
with
\beqano
\ovl f^{\rm av}=\ovl C_0(\L)+\O\cdot\t+\frac{1}{2}\t\cdot{\rm T}(\L)\t+\sum_{j=3}^s\ovl\cP(\t;\L)+\ovl\cR_{2s+2}(\L,\ovl\eta,\ovl\xi,\ovl p, \ovl q)
\eeqano
with 
$$\t_1=\ovl\t_1=\frac{\ovl\eta_1^2+\ovl\xi_1^2}{2}\ ,\quad \t_2=-\ovl\t_2=\frac{\ovl\eta_2^2+\ovl\xi_2^2}{2}\ ,\quad \t_3=-\ovl\t_3=\frac{\ovl p^2+\ovl q^2}{2}$$
and $\O$, {\rm T} as in~\equ{TT}. The proof that $\ovl\phi$ is real--analytic {in $\ovl\cM_{\ovl\varepsilon}$} is standard and therefore is omitted.\\
(iii) follows from (i) and (ii).
(iv) is analogous to Refs.\cite{chierchiaPi11b, chierchiaPi11c}.


\subsubsection{\label{sec: full Birkhoff}Birkhoff theory {the} fully reduced retrograde problem}

{The construction of the Birkhoff normal form for the averaged, fully reduced system in the case of the prograde problem in has been discussed in Ref.\cite[Sections 9 and 10]{chierchiaPi11b}}. We recall the basic steps:
\begin{itemize}
\item[$\star$] Passage to a set of canonical coordinates on the manifold with constant $\GG=|\CC|$. This is accomplished via a transformation $\hat\phi$ described in Ref.\cite[Section 9]{chierchiaPi11b} which conjugates a ``Birkhoff--normalized partially reduced system'' $\breve\cH$ discussed in Ref.\cite[Equations (7.28)--(7.30)]{chierchiaPi11b} to a fully reduced system $\hat\cH_\GG$, as described in  Ref.\cite[Remark 9.1]{chierchiaPi11b};
\item[$\star$] Construction, via an Implicit Function Theorem argument (e.g., Ref.\cite[Proposition 10.1]{chierchiaPi11b}) of a canonical transformation $\check\phi$ which conjugates $\hat\cH_\GG$ to $\check\cH_\GG$, where the average $\check\cH_\GG$ with respect the the angles $\check\l$ appearing in Ref.\cite[Equation (10.10)]{chierchiaPi11b} is in Birkhoff normal form of order 4.
\end{itemize}
 In the case of the retrograde problem,  one can mimic  the procedure above, or, equivalently,  exploit the results of Ref.\cite{chierchiaPi11b} via Proposition~\ref{signs1}. Namely, if $\hat\phi^\complex:\ (\L,\GG, \hat\l, \hat{\rm g}, \hat t, \hat t^*)\to (\L,\breve\l, \breve t, \breve t^*)$, $\check\phi^\complex:\ (\L,\GG, \check\l, \check{\rm g}, \check t, \check t^*)\to (\L,\GG, \hat\l, \hat{\rm g}, \hat t, \hat t^*)$ denote the transformations induced by the transformations
 $\check\phi:\ (\L,\GG, \hat\l, \hat{\rm g}, \hat \eta, \hat\xi)\to (\L,\breve\l, \breve \eta, \breve \xi, \breve p, \breve q)$, $\check\phi:\ (\L,\GG,  \check\l, \check{\rm g}, \check \eta, \check \xi)\to (\L,\GG, \hat\l, \hat{\rm g}, \hat \eta, \hat\xi)$ of Ref.\cite{chierchiaPi11b} on the  complex coordinates 
 related to the respective real coordinates via   relations  to \equ{PROld+};  $\phi_2^-$ is as in \equ{involution}; $\dot\phi^\complex:=\phi_2^-\circ\hat\phi^\complex\circ\check\phi^\complex\circ\phi_2^-:\ (\L,\dot\l, \dot t, \dot t^*)\to (\L, \ovl\l, \ovl t, \ovl t^*)$ then the transformation $\dot\phi:\ (\L,\GG, \dot\l, \dot\g, \dot\eta, \dot\xi)\to(\L,\ovl\l, \ovl\eta, \ovl\xi, \ovl p, \ovl q)$  induced by $\dot\phi^\complex$, such in a way that 
 $(\L,\GG, \dot\l, \dot\g, \dot\eta, \dot\xi)$, $(\L,\ovl\l, \ovl\eta, \ovl\xi, \ovl p, \ovl q)$ are related to 
$ (\L,\dot\l, \dot t, \dot t^*)$, $(\L, \ovl\l, \ovl t, \ovl t^*)$, respectively, via relations \equ{PR1},
turns to be real on a real domain and
 makes the following proposition true.

 \begin{proposition}[BNF for the fully reduced retrograde problem]\label{BNF}
Fix $4\le s\in \natural$.
There exist $0<\dot\varepsilon_1<\frac{\dot\varepsilon_2}{4}$ such that, if
\beqa{LM}
\dot\cL_s&:=&\Big\{(\L, \GG)=(\L_1,\L_2, \GG)\in \cL_s\times \real_+:\ \L_1>\L_2\ ,\quad \GG<\L_1-\L_2\nonumber\\
&& \varrho(\L,\GG):=\sqrt{2(\L_1-\L_2-\GG)}\in \Big(4\dot\varepsilon_1, \dot\varepsilon_2\Big)\Big\}\nonumber\\
&&\dot\cM_{\dot\varepsilon}:=\dot\cL_s\times \torus^3\times B^{4}_{\dot\varepsilon}(0)
\eeqa
one can find 
a real--analytic and symplectic transformation 
\beqa{fully reduced retrograde}
\dot\phi:\quad  \dot\cM_{\dot\varepsilon_1}=\dot\cL_s\times \torus^3\times B^{4}_{\dot\varepsilon_1}(0)&\to& \ovl\cM_{\ovl\varepsilon}:=\cL_{}\times \torus^2\times B^{6}_{\ovl\varepsilon}(0)\nonumber\\
(\L,\GG, \dot\l, \dot\g,  \dot\eta, \dot\xi)&\to&(\L,\dot\l+\dot\f(\L,\dot\eta, \dot\xi), \dot B(\dot\eta, \dot\xi;\L))
\eeqa such that \beqa{eq: normalized orc}\dot\HH_\GG(\L,\dot\l,\dot\eta, \dot\xi)=\ovl\HH\circ\dot\phi(\L,\dot\l,\dot\eta, \dot\xi;\GG)=h_{\rm k}(\L)+\m\dot f(\L,\dot\l,\dot\eta, \dot\xi;\GG)\eeqa
is $\dot\g$--independent, and
 $\dot f^{\rm av}(\L,\dot \eta,\dot\xi;\GG):=\frac{1}{4\p^2}\int_{\torus^2}\dot f(\L,\dot\l,\dot\eta, \dot\xi;\GG)d\dot\l$
takes the form
\beqa{dot T}
\dot f^{\rm av}(\L,\dot \eta,\dot\xi;\GG)&=&\dot C(\L;\GG)+\dot\O(\L;\GG)\cdot\dot\t+\frac{1}{2}\, \dot\t \dot{\rm T}(\L;\GG)\cdot\dot\t+\sum_{j=3}^{s}\dot\cP_j(\dot\t;\L,\GG)\nonumber\\
&+&
{\rm O}_{2s+1}(\L,\dot\eta,\dot\xi;\GG)
\eeqa
where $\dot\t_i:=\frac{\dot\eta_i^2+\dot\xi_i^2}{2}$, and
\beqa{reduced invariants orc}
&&\dot\O_i(\L_1,\L_2)=s_i\hat\O_i(\L_1,-\L_2)\nonumber\\
&&\dot{\rm T}_{ij}(\L_1,\L_2;\GG)=s_is_j\hat{\rm T}_{ij}(\L_1,-\L_2)+{\rm O}(\dot\varepsilon_2^2)\ ,\ 1\le i\le j\le 2\eeqa
with $s_1=-s_2=1$
and
\beqa{reduced invariants+}
&&\hat\O_i(\L_1,\L_2)=\breve\O_i(\L_1,\L_2)-\breve\O_3(\L_1,\L_2)\nonumber\\
&&\hat{\rm T}_{ij}(\L_1,\L_2):=\breve{\rm T}_{ij}(\L_1,\L_2)-\breve{\rm T}_{i3}(\L_1,\L_2)-\breve{\rm T}_{3j}(\L_1,\L_2)+\breve{\rm T}_{33}(\L_1,\L_2)\ .\eeqa
\end{proposition}

{\begin{remark}\rm
The analytical expression $\dot\O_1(\L_1,\L_2)$, $\dot\O_2(\L_1,\L_2)$ which is available accordingly to \equ{reduced invariants+} and Ref.\cite[Equations (7.1), (7.5) with $n=2$]{chierchiaPi11b} 
$$\dot\O_1=-\frac{3}{4}m_1 m_2\frac{a_1^2}{a_2^3}(\frac{2}{\L_1}-\frac{1}{\L_2})(1+{\rm O}\left(\frac{a_1^2}{a_2^2}\right))\quad \dot\O_2=\frac{3}{4}m_1 m_2\frac{a_1^2}{a_2^3}(\frac{1}{\L_1}-\frac{2}{\L_2})(1+{\rm O}\left(\frac{a_1^2}{a_2^2}\right))$$
shows that the $\dot\O_i$'s
do not satisfy, identically, linear combinations up to any prefixed order in $\dot\cL_{s,\GG}$. This is not used in the paper, but is useful
to be known if one wants to obtain a stronger  result concerning the measure of the invariant set (Ref.\cite[Theorem 1.2]{chierchiaPi10})
\end{remark}}
\subsubsection{\label{sec: torsion} Torsion}
{ Define, for a fixed $\GG\in \real$, \beqa{LsG}\cL_{s,\GG}:=\big\{(\L_1,\L_2):\ (\L_1,\L_2, \GG)\in \dot\cL_s\big\}\eeqa
where $\dot\cL_s$ is as in \equ{Ls}.
In this section, we first compute the  the matrix $\dot{\rm T}$ in \equ{reduced invariants orc} and next we check its non--singularity 
on  $\dot\cL_{s,\GG}$. We proceed in three steps.}

\vskip.1in
\paragraph{Computation of $\breve{\rm T}$}
Note that we cannot use the asymptotics of the matrix $\breve{\rm T}$ given in Ref.\cite[Equation (8.6)]{chierchiaPi11b}, because that was obained for $\L_1\ll\L_2$, while, for the purposes of the paper, we need $\L_1>\L_2$.\\
By~\equ{TT}, we evaluate the $\breve{\rm T}_{ij}(\L_1,\L_2)$'s  first. 
According to~Refs.\cite{pinzari-th09, chierchiaPi11b}, they  have the form (neglecting to write the arguments)
 \beqa{three bodies}
\breve{\rm T}_{11}&=&\frac{4m_1m_2}{(1+d^2)^2}\Big[\frac{{\rm r}_1(a_1,a_2)}{\L_1^2}+\frac{d^4{\rm r}_1(a_2,a_1)}{\L_2^2}+\frac{2d^2{\rm r}_2(a_1,a_2)}{\L_1\L_2}\\
&&\phantom{AAAAAA} -\frac{2d{\rm r}_3(a_1,a_2)}{\L_1\sqrt{\L_1\L_2}}-\frac{2d^3{\rm r}_3(a_2,a_1)}{\L_2\sqrt{\L_1\L_2}}+\frac{d^2{\rm r}_4(a_1,a_2)}{\L_1\L_2}\Big]\nonumber\\
\breve{\rm T}_{12}&=&\frac{4m_1m_2}{(1+d^2)^2}\Big[\frac{(1-d^2)^2{\rm r}_2(a_1,a_2)}{\L_1\L_2}+\frac{2d^2{\rm r}_1(a_1,a_2)}{\L_1^2}+\frac{2d^2{\rm r}_1(a_2,a_1)}{\L_2^2}\nonumber\\
&&\phantom{AAAAAA}+\,\frac{2d(1-d^2){\rm r}_3(a_1,a_2)}{\L_1\sqrt{\L_1\L_2}}-\frac{2d(1-d^2){\rm r}_3(a_2,a_1)}{\L_2\sqrt{\L_1\L_2}}+\frac{(1-6d^2+d^4){\rm r}_4(a_1,a_2)}{4\L_1\L_2}\Big]\nonumber\\
\breve{\rm T}_{22}&=&\frac{4m_1m_2}{(1+d^2)^2}\Big[\frac{{\rm r}_1(a_2,a_1)}{\L_2^2}+\frac{d^4{\rm r}_1(a_1,a_2)}{\L_1^2}+\frac{2d^2{\rm r}_2(a_1,a_2)}{\L_1\L_2}\nonumber\\
&&\phantom{AAAAAA}-\frac{2d{\rm r}_3(a_2,a_1)}{\L_2\sqrt{\L_1\L_2}}-\frac{2d^3{\rm r}_3(a_1,a_2)}{\L_1\sqrt{\L_1\L_2}}+\frac{d^2{\rm r}_4(a_1,a_2)}{\L_1\L_2}\Big]\nonumber
\\
\breve{\rm T}_{13}&=&\frac{m_1m_2}{1+d^2}\Big[\frac{1}{\L_1}\Big(\frac{1}{\L_1}+\frac{1}{\L_2}\Big)({\rm s}_1(a_1,a_2)+{\rm s}_1^*(a_1,a_2))-\Big(\frac{1}{\L_1^2}+\frac{d^2}{\L_2^2}\Big)C_1(a_1,a_2)\nonumber\\
&&\phantom{AAAAAA}-\frac{d}{\sqrt{\L_1\L_2}}\Big(\frac{1}{\L_1}+\frac{1}{\L_2}\Big)({\rm s}_2(a_1,a_2)+{\rm s}_2^*(a_1,a_2))\nonumber\\
&&\phantom{AAAAAA}+\,\frac{d^2}{\L_2}\Big(\frac{1}{\L_1}+\frac{1}{\L_2}\Big)({\rm s}_1(a_2,a_1)+{\rm s}_1^*(a_2,a_1))\Big]\nonumber\\
\breve{\rm T}_{23}&=&\frac{m_1m_2}{1+d^2}\Big[\frac{1}{\L_2}\Big(\frac{1}{\L_1}+\frac{1}{\L_2}\Big)({\rm s}_1(a_2,a_1)+{\rm s}_1^*(a_2,a_1))-\Big(\frac{1}{\L_2^2}+\frac{d^2}{\L_1^2}\Big)C_1(a_1,a_2)\nonumber\\
&&\phantom{AAAAAA}+\,\frac{d}{\sqrt{\L_1\L_2}}\Big(\frac{1}{\L_1}+\frac{1}{\L_2}\Big)({\rm s}_2(a_1,a_2)+{\rm s}_2^*(a_1,a_2))\nonumber\\
&&\phantom{AAAAAA}+\,\frac{d^2}{\L_1}\Big(\frac{1}{\L_1}+\frac{1}{\L_2}\Big)({\rm s}_1(a_1,a_2)+{\rm s}_1^*(a_1,a_2))\Big]\nonumber\\
\breve{\rm T}_{33}&=&4m_1m_2\Big[{\rm r}_1^*(a_1,a_2)\Big(\frac{1}{\L_1}+\frac{1}{\L_2}\Big)^2+\frac{C_1(a_1,a_2)}{4\L_1\L_2}\Big]\nonumber
\eeqa
%
where 
$d$ depends on $\L_1$, $\L_2$, $a_1$, $a_2$, while 
$C_1$,  ${\rm r}_1$,  ${\rm r}_2$ ${\rm r}_3$, ${\rm r}_4$, ${\rm r}_1^*$, ${\rm r}_2^*$, , ${\rm s}_1$, ${\rm s}_2$,  ${\rm s}_1^*$, ${\rm s}_2^*$ are
 functions of $a_1$, $a_2$ only (in turn related to $\L_1$, $\L_2$ via~\equ{aL}), expressed in terms of
the Laplace coefficients. Moreover, it follows from~Ref.\cite[Proof of Proosition 5.2]{pinzari-th09} that
\beqa{dnew}d=d(\frac{\sqrt{\L_1\L_2}}{|\L_1-\L_2|}b(\a))\eeqa
with
 $$d(x)=\frac{x}{\sqrt{1+x^2}+1}\ ,\quad b(\a)=4\frac{b_{3/2}^\ppd(\a)}{b_{3/2}^\ppu(\a)}={\rm O}(\a)\ .$$
 We recall that $d$ is defined so that the matrix $\widehat U$ in~\equ{hat U} is
 $$
 \left(
 \begin{array}{cccc}
 1&-\frac{d}{\sqrt{1+d^2}}\\
 \frac{d}{\sqrt{1+d^2}}&1
 \end{array}
 \right)\ .$$
 As in ~Refs.\cite{pinzari-th09,chierchiaPi11b}, we shall use the expansions of the functions above in terms of $\a$
\beqa{neglect}
\begin{array}{lllll}
 {\rm r}_1(a_1,a_2)=\frac{3}{16a_2}\Big(\a^2+{\rm O}\Big(\a^4\Big)\Big)\qquad  &{\rm r}_1(a_2,a_1)=-\frac{3}{4a_2}\Big(\a^2+{\rm O}\Big(\a^4\Big)\Big)\\ 
{\rm r}_2(a_1,a_2)={\rm r}_2(a_2,a_1)=-\frac{9}{16a_2}\Big(\a^2+{\rm O}\Big(\a^4\Big)\Big) &{\rm r}_3(a_1,a_2)={\rm O}\Big(\frac{\a^3}{a_2}\Big)={\rm r}_3(a_2,a_1)\\ 
{\rm r}_4(a_1,a_2)={\rm O}\Big(\frac{\a^4}{a_2}\Big) &C_1(a_1,a_2)=-\frac{3}{4a_2}\Big(\a^2+{\rm O}\Big(\a^4\Big)\Big)\\
{\rm r}_1^*(a_1,a_2)=-\frac{3}{16a_2}\Big(\a^2+{\rm O}\Big(\a^4\Big)\Big)\quad &{\rm s}_1(a_1,a_2)=\frac{3}{a_2}\left(\a^2+{\rm O}\left(\a^4\right)\right)\\ 
{\rm s}_1(a_2,a_1)=\frac{9}{8a_2}\left(\a^2+{\rm O}\left(\a^4\right)\right)& {\rm s}^*_1(a_1,a_2)=-\frac{3}{4a_2}\left(\a^2+{\rm O}\left(\a^4\right)\right)\\ 
{\rm s}^*_1(a_2,a_1)=\frac{9}{8a_2}\left(\a^2+{\rm O}\left(\a^4\right)\right)& {\rm s}_2(a_1,a_2),\ {\rm s}_2^*(a_1,a_2)={\rm O}\left(\frac{\a^3}{a_2}\right)\\ 
\end{array}
\eeqa
Moreover, by  the definition of $d$ in~\equ{dnew}, $d={\rm O}(\a\sqrt{t})$, with $t:=\frac{\L_2}{\L_1}$. We thus obtain
\beqa{old torsion}
\breve{\rm T}={m_1m_2}\frac{\a^2}{a_2\L_2^2}\left(
\begin{array}{ccccc}
\frac{3}{4}t^2&\quad&-\frac{9}{4}t&\quad&3t^2-\frac{9}{4}t\\
-\frac{9}{4}t&&-3&&\frac{9}{4}t+3\\
3t^2-\frac{9}{4}t&&\frac{9}{4}t+3&&-\frac{3}{4}(1+t)^2-\frac{3}{4}t\\
\end{array}
\right)(1+{\rm O}(\a^2))
\eeqa

\vskip.1in
\paragraph{Computation of $\hat{\rm T}(\L_1,\L_2)$}
Using~\equ{old torsion} and ~\equ{reduced invariants+}, we obtain
\beqano
\hat{\rm T}(\L_1,\L_2)=-\frac{3}{4}{ m_1 m_2}\frac{\a^2}{a_2\L_2^2}\left(
\begin{array}{llll}
1-3t+8 t^2&\quad & 5+6 t+5 t^2\\
5+6 t+5 t^2&&13+9t+t^2\\\end{array}
\right)(1+{\rm O}(\a^2))\ .
\eeqano
\vskip.1in
\paragraph{Computation of $\dot{\rm T}(\L_1,\L_2;\GG)$ and check of its non--singularity}
The matrix $\dot{\rm T}(\L_1,\L_2;\GG)$ defined in \equ{reduced invariants orc} is given by
\beqano
\dot{\rm T}(\L_1,\L_2;\GG)=-\frac{3}{4}{ m_1 m_2}\frac{\a^2}{a_2\L_2^2}\left(
\begin{array}{llll}
1+3t+8 t^2&\quad & -5+6 t-5 t^2\\
-5+6 t-5 t^2&&13-9t+t^2\\\end{array}
\right)(1+{\rm O}(\a^2))+{\rm O}(\varepsilon_2^2)\ .
\eeqano
It may vanish at most on a finite number of sub--manifolds of $\dot\cL_\GG$:
$$\det\dot{\rm T}(\L_1,\L_2;\GG)=-\frac{9}{16}\Big({ m_1 m_2}\frac{\a^2}{a_2\L_2^2}\Big)^2\big(
12-90 t+8 t^2+9 t^3+17 t^4
\big)(1+{\rm O}(\a^2))+{\rm O}(\varepsilon_2^2)\ .$$

\subsection{\label{domain of rps+}Existence of quasi--periodic motions}

Let us now take, as mentioned, $\a_+<\a_s$, where $\a_s$ is as in~\equ{Ls}. Let  $\dot\HH$  be as in~\equ{eq: normalized orc}. We aim to apply the following result to $\dot\HH$. Note that assumption (A${}_1$) is trivially satisfied by
$\hh_{\rm k}$, while (A${}_2$) and (A${}_3$) have been discussed in Sections~\ref{sec: full Birkhoff} and~\ref{sec: torsion}.

\begin{theorem}[V.I.Arnold,~Refs.\cite{arnold63,chierchiaPi10}]\label{thm:*simplified}
Let ${\cal P}_\varepsilon:=V\times\torus^{n_1}\times B_{\varepsilon}^{2n_2}$,
where $V$ is an open, bounded, connected set of $\real^{n_1}$ and $B^{2n_2}_{\varepsilon}$ is a $2n_2$--dimensional ball of radius $\varepsilon$ centered at the origin. Let $\varepsilon_0>0$ and let $H(I,\varphi,p,q;\m)=H_0(I)+\m\,P(I,\varphi,p,q;\m)$ be a real--analytic Hamiltonian on $\cP_{\varepsilon_0}$, endowed with  the standard symplectic form $dI\wedge d\varphi+dp\wedge dq$.
Assume that $H$ verifies the following non--degeneracy assumptions: 
\begin{itemize}
\item[\rm (A${}_1$)]  
$I\in V\to \partial_I H_0$ is a diffeomorphism; 

\item[\rm (A${}_2$)]  
$ P_{\!\rm av} (p,q;I)=P_0(I) +\sum_{i=1}^{n_2}\O_i(I)r_i+\frac{1}{2}\sum_{i,j=1}^{n_2} \b_{ij}(I)r_i r_j+{\rm O}_{5}(p,q;I)$,  with  $r_i:=\frac{p_i^2+q_i^2}{2}$; 

\item[\rm (A${}_3$)]  
 $|\det \b(I)|\ge \const>0$ for all $I\in V$.
\end{itemize}
Then, there exist positive numbers $\varepsilon_*<\varepsilon_0$, $C_*$ and $c_*$ such that, for 
\beqa{epsmu}0<\varepsilon<\varepsilon_*\ , \qquad 
0<\m<\frac{\varepsilon^6}{ (\log \varepsilon^{-1})^{c_*}}\ ,
\eeqa
one can find a 
set $\cK\subset  \cP_\varepsilon$  formed by the union of $H$--invariant $n$--dimensional Lagrangian  tori, on which the $H$--motion is analytically conjugated to  linear Diophantine quasi--periodic motions with frequencies $(\o_1,\o_2)\in\real^{n_1}\times\real^{n_2}$ with $\o_1=O(1)$ and $\o_2=O({\m})$. The  set $\cK$  has positive Liouville--Lebesgue measure and satisfies
\beqa{meas est}
\meas \cP_\varepsilon>\meas \cK> \Big(1- C_* \sqrt{\varepsilon}\Big) \meas \cP_\varepsilon\ .
\eeqa
\end{theorem}

We now aim to apply Theorem \ref{thm:*simplified} to the Hamiltonian $\dot\HH_\GG$ in \equ{eq: normalized orc}, for any fixed $\GG\in \real_+$. So,
if $\dot\varepsilon_1$  as in Proposition~\ref{BNF}, 
$\varepsilon\in (0,\dot\varepsilon_1)$;   $\dot\cL_{s,\GG}$ is as in \equ{LsG} and
$\dot\cM_{\GG, \varepsilon}:=\dot\cL_{s,\GG}\times \torus^2\times B^4_\varepsilon
$, we obtain the following
\begin{theorem}\label{stable tori} 
{
Let $\varepsilon$ and $\m$ verify \equ{epsmu}.
Then the set $\dot\cM_{\GG, \varepsilon}$ contains a positive measure set $\dot\cK_\GG$ which is $\dot\HH_\GG$ invariant and is formed by the union of $4$--dimensional Lagrangian, real--analytic tori on which the $\dot\HH_\GG$ motions are analytically conjugated to linear quasi--periodic motions. Furthermore, $\dot\cK_\GG$ satisfies the bound in \equ{meas est}, with $\cP_\varepsilon$ replaced by $\dot\cM_{\GG, \varepsilon}$.
}
\end{theorem}

\section{\label{generalities}The reduction of perihelia}
This section  contains a new coordinate system, which, just like the {\rm jrd} coordinates in~\equ{coordinates}
 is of the form ``action--angle'' and reduces completely the number of degrees of freedom, and, in addition, enjoys nice and useful parity properties.
It is  available for any number of planets (Ref.\cite{pinzari15}), but we discuss it here only for two. 

\vskip.2in
We  call {\it Perihelia reduction for the three-body problem} 
the set of coordinates $${\rm p}:=(\L_1,\L_2,\GG_2,\Theta,\GG,\ZZ, \ell_1,\ell_2, {\rm g}_2, \vartheta, {\rm g}, \zeta)$$
 defined as follows.\\
Consider the sets of $(x,y)$ such that the  Keplerian motions  $t\to (x^\ppj(t), y^\ppj(t))$ 
generated by the  Hamiltonians~\equ{K***} starting from $(x^\ppj,y^\ppj)$ are ellipses with non--vanishing eccentricities. Let $a_i$ denote their semi--major axes; $\CC^\ppi:=x^\ppi\times y^\ppi$ the angular momenta; $\PP^\ppi$, with $|\PP^\ppi|=1$, the direction of the $i^{\rm th}$ perihelion. 
Define the ``${\rm p}$--nodes''
\beqa{good nodes}
 \n_1:=k^\ppt\times {\rm C},\qquad{\rm n}_1:={\rm C}\times {\rm P}^\ppu,\qquad   \n_2:={\rm P}^{(1)}\times {\rm C}^\ppd, \quad {\rm n}_2={\rm C}^\ppd\times {\rm P}^\ppd
\eeqa
and assume that they do not vanish. If, as in the definition of ${\rm jrd}$, for any three vectors $u$, $v$, $w\in \real^3$, with $u$, $v\perp$ $w$, let $\a_w(u,v)$ denote the oriented angle formed from  $u$ to  $v$, relative to  the positive direction established by $w$, then  define
\beqa{belle*}
\begin{array}{llllrrr}
 \left\{
\begin{array}{lrrr}
 \ZZ:={\rm C}\cdot k^\ppt \\
 \Theta:={\rm C}\cdot {\rm P}^{(1)}={\rm C}^\ppd\cdot {\rm P}^{(1)} \\
{\rm G}:=\|{\rm C}\| \\
{\rm G}_2:=\|{\rm C}^\ppd\|\\
\L_j:={\rm  M}_j\sqrt{{\rm  m}_j a_j}
\end{array}
\right.\qquad
\left\{
\begin{array}{lrrr}
\zeta:=\a_{k^\ppt}(k^\ppu, \n_1)\qquad& \\
 \vartheta:=\a_{{\rm P}^{(1)}}({\rm n}_1, \n_2)&\\
 {\rm g}:=\a_{{\rm C}}(\n_1, {\rm n}_1)\qquad&\\
 {\rm g}_2:=\a_{{\rm C}^\ppd}(\n_2, {\rm n}_2)&\\
\ell_j := {\rm mean\ anomaly\ of}\ x^\ppj& 
\end{array}
\right.
\end{array}
\eeqa
with $j=1, 2$.
Note that such definitions, together with  assumptions of non--vanishing of the nodes~\equ{good nodes}, imply 
that
\beqa{p conditions}
\min\Big\{\GG_2,\ \GG\Big\}>0\ ,\quad -\GG<\ZZ<\GG\ ,\quad -\min\Big\{\GG_2,\ \GG\Big\}<\Theta<\min\Big\{\GG_2,\ \GG\Big\}\ .
\eeqa
 In Section~\ref{appendix D} we shall show that the map 
\beqa{p}\phi_{{\rm p}}: \qquad {\rm p}:=(\ZZ,\Theta,\chi,\L,\zeta,\vartheta,\k,\ell)\in \real^6\times \torus^6\to\cC_{\rm art} (y_{\rm p}^\ppu, y_{\rm p}^\ppd,x_{\rm p}^\ppu,x_{\rm p}^\ppd)\in \real^{12}\eeqa
preserves that standard 2--form:
$$d\ZZ\wedge d\zeta+d \GG\wedge d{\rm g}+d\Theta\wedge d\vartheta+d\GG_2\wedge d{\rm g_2}+d\L_1\wedge d\ell_1+d\L_2\wedge d\ell_2=\sum_{j=1}^2\sum_{i=1}^3dy^\ppj_{{\rm p},i} \wedge d y^\ppj_{{\rm p},i}\ .$$
\subsubsection{\label{phiP}Analytical expression of $\phi_{\rm p}$}
The explicit form of $\phi_{\rm p}$ is as follows. Let  ${\rm e}_2$, $\zeta_2$, $i$,  $R_1$, $R_3$ be as in~\equ{aL},~\equ{ej},~\equ{zetaj},~\equ{D incli},~\equ{R1R3}. Let $\iota_1$, $\iota_2$ be the convex angles
\beqano \cos \iota_1=\frac{\Theta}{\GG}\ ,\quad \cos \iota_2=\frac{\Theta}{\GG_2}\ .\eeqano
By~\equ{p conditions}, one has $\iota_1$, $\iota_2\in (0,\p)$.
Then the expression of $\CC$, $\CC^\ppd$, $\PP^\ppu$, $\PP^\ppd$ in terms of ${\rm p}$ are
\beqa{P formulae}
&&\CC=\GG R_3(\zeta)R_1(i)e_3\ ,\quad \PP^\ppu=R_3(\zeta)R_1(i) R_3({\rm g})R_1(\iota_1)e_3\nonumber\\
&& \CC^\ppd=\GG_2 R_3(\zeta)R_1(i) R_3({\rm g})R_1(\iota_1)R_3(\vartheta)R_1(\iota_2)e_3\nonumber\\
&& \PP^\ppd=R_3(\zeta)R_1(i) R_3({\rm g})R_1(\iota_1)R_3(\vartheta)R_1(\iota_2) R_3({\rm g}_2-\p/2)e_1
\eeqa
while $\CC^\ppu$ is found via $$\CC^\ppu=\CC-\CC^\ppd=R_3(\zeta)R_1(i)\Big(\GG {\bf 1}-\GG_2 R_3({\rm g})R_1(\iota_1)R_3(\vartheta)R_1(\iota_2)\Big) e_3\ . $$
In particular, $\GG_1:=\|{\rm C}^\ppu\|$ is not an action coordinate, but has the expression
\beqa{C}
\GG_1=\sqrt{\GG^2+\GG_2^2-2\Theta^2+2\sqrt{\GG^2-\Theta^2}\sqrt{\GG_2^2-\Theta^2}\cos\vartheta}\ .
\eeqa
Such expression allows to find ${\rm e}_1$ via~\equ{ej} as
$${\rm e}_1=\sqrt{1-\frac{\GG^2+\GG_2^2-2\Theta^2+2\sqrt{\GG^2-\Theta^2}\sqrt{\GG_2^2-\Theta^2}\cos\vartheta}{\L_1^2}}$$
Letting now
\beqa{Q1Q2}
{\rm Q}^\ppu:=\frac{\CC^\ppu}{\GG_1}\times \PP^\ppu\ ,\quad 
{\rm Q}^\ppd:=\frac{\CC^\ppd}{\GG_2}\times \PP^\ppd
\eeqa
then, $x^\ppj_{\rm p}:=x^\ppj$, $y^\ppj_{\rm p}:=y^\ppj$ are, classically, given by
\beqa{xyP}
\arr{x^\ppj=\frac{\L^2_j}{\MM_j\mm_j^2}\Big((\cos\zeta_j-{\rm e}_j)\PP^\ppj+\sqrt{1-{\rm e}_j^2}\sin\zeta_j{\rm Q}^\ppj\Big)\\
y^\ppj=\frac{\MM_j\mm_j^2}{\L_j}\Big(-\frac{\cos\zeta_j}{1-{\rm e}_j\cos\zeta_j}\PP^\ppj+\sqrt{1-{\rm e}_j^2}\frac{\cos\zeta_j}{1-{\rm e}_j\cos\zeta_j}{\rm Q}^\ppj\Big)
}
\eeqa
where $\zeta_j$ solves~\equ{zetaj}.
Observe that, while the expression of ${\rm Q}^\ppd$in~\equ{Q1Q2} is relatively simple:
\beqano
{\rm Q}^\ppd&=&\frac{\CC^\ppd}{\GG_2}\times \PP^\ppd=R_3(\zeta)R_1(i) R_3({\rm g})R_1(\iota_1)R_3(\vartheta)R_1(\iota_2) R_3({\rm g}_2)e_1
\eeqano
the one of ${\rm Q}^\ppu$ is much more involved. Fortunately, it will be not needed in the paper.
\subsubsection{Singularities}
By definition, the singularities of $\phi_{\rm p}$ arise when
\beqano
&&\ZZ=\GG\quad {\rm or}\quad \ZZ=\GG\quad {\rm or}\quad \Theta=\GG\quad {\rm or}\quad\Theta=\GG_2\quad {\rm or}\quad\GG_2=\L_2\nonumber\\
&&\GG_1=\sqrt{\GG^2+\GG_2^2-2\Theta^2+2\sqrt{\GG^2-\Theta^2}\sqrt{\GG_2^2-\Theta^2}\cos\vartheta}=\L_1\quad {\rm or}\quad \GG_2=0\quad {\rm or}\quad\GG=0\nonumber\\
&&{\rm or}\quad (\GG_2=\GG\quad {\rm and}\quad \vartheta=\p)\ .
\eeqano
Note that the last line corresponds to the vanishing of $\GG_1$ in~\equ{C}.
These formulae show that, compared to $\phi_{\rm jrd}$, $\phi_{\rm p}$ has the advantage that {\it  there are not singularities for vanishing inclinations}.
As a counterpart, while the singularity for $\GG_2=\L_2$ can be eliminated switching to the regularized (planar) Poincar\'e coordinates
\beqano
(\L_2,\ell_2,\GG_2,{\rm g}_2)\to (\L_2,{\rm l}_2, u_2,v_2):=(\L_2,\ell_2+{\rm g}_2, \sqrt{2(\L_2-\GG_2)}\cos{\rm g}_2, -\sqrt{2(\L_2-\GG_2)}\sin{\rm g}_2)\eeqano
for the spatial problem,
a regularization for the singularity for  $\GG_1=\L_1$ does not seem to be possible, due to the non--linear expression of $\GG_1$ in~\equ{C}.

\subsubsection{Symmetries and planar equilibria}
Let
\beqa{p Ham}
{\rm H}_{{\rm p}}&:=&\HH\circ\phi_{\rm p}=\hh_{\rm k}(\L_1,\L_2)+\m f_{\rm p}(\L_1,\L_2,\GG_2,\Theta;\ell_1,\ell_2,{\rm g}_2,\vartheta; \GG)\eeqa
the Hamiltonian $\HH$ expressed in terms of ${\rm p}$. As well as $\HH_{\rm jrd}$ in~\equ{jrd Ham},  $\HH_{\rm p}$ is independent of $\GG$, $\ZZ$ and ${\rm g}$. However, as an advantage with respect to the former, the symmetries by reflections induce a useful parity property to $\HH_{\rm p}$, which does not exist for $\HH_{\rm jrd}$. Indeed, it turns out that the transformation
\beqano
\cR_2^-:\quad (\L_1,\L_2, \GG_2,\Theta, \GG,\ZZ, \ell_1,\ell_2, {\rm g}_2, \vartheta, {\rm g}, \zeta)\to (\L_1,\L_2, \GG_2,-\Theta, \GG,-\ZZ, \ell_1,\ell_2, {\rm g}_2, -\vartheta, {\rm g}, -\zeta)\eeqano
corresponds, in the Cartesian coordinates in~\equ{helio}, to
\beqano
\cR_2^-:\quad(y^\ppj_1,y^\ppj_2, y^\ppj_3, x^\ppj_1,x^\ppj_2, x^\ppj_3)\to(y^\ppj_1,-y^\ppj_2, y^\ppj_3, x^\ppj_1,-x^\ppj_2, x^\ppj_3)\qquad j=1, 2\eeqano
This can be easily verified using the formulae given in Section~\ref{phiP}. Then, the $f_{\rm p}$
enjoys the following parity property
\beqano
f_{\rm p}(\L_1,\L_2,\GG_2,\Theta;\ell_1,\ell_2,{\rm g}_2,\vartheta; \GG)=f_{\rm p}(\L_1,\L_2,\GG_2,-\Theta;\ell_1,\ell_2,{\rm g}_2,-\vartheta; \GG)\quad {\rm mod}\ 2\p\ .\eeqano
This equalities implies that the three  manifolds
\beqa{equilibria}
\begin{array}{llll}
(\uparrow\,\uparrow)&:=\Big\{{\rm p}: \ (\Theta, \vartheta)=(0,\p)\quad{\rm and}\quad \GG>\GG_2\Big\}; \\
(\downarrow\,\uparrow)&:=\Big\{{\rm p}: \ (\Theta, \vartheta)=(0,\p)\quad {\rm and}\quad \GG<\GG_2\Big\}; \\
(\uparrow\,\downarrow)&:= \Big\{{\rm p}: \ (\Theta, \vartheta)=(0,0)\Big\}.
\end{array}
\eeqa
are invariant to $f_{\rm p}$. Such manifolds correspond to planar motions, with the directions of the left, right arrows denoting the direction of the angular momentum of the inner, outer planet, respectively. The manifold  $(\uparrow\,\downarrow)$ corresponds to planar {\sc orc} motions, which the object of study of this paper.

\subsubsection{\label{domain}A domain of regularity including ${(\uparrow\downarrow)}$}

In this section we establish  a suitable domain  which includes the invariant manifold $(\uparrow\downarrow)$ where ${\rm H}_{{\rm p}}$ is regular.

We  check below that the following domain is suited to the scope: 
\beqano
\cD(\GG)&:=&\Big\{(\L_1, \L_2, \GG_2,\Theta,  \vartheta): \ (\L_1, \L_2, \GG_2)\in \cA(\GG), (\Theta,\vartheta)\in \cB(\GG_2,\GG)\Big\}\times\torus^3.
\eeqano
where, if $\cL$ is as in~\equ{L0},  
\beqa{B}
\cA(\GG)&:=&\Big\{(\L_1, \L_2, \GG_2): \ (\L_1,\L_2)\in \cL(\GG), \GG_2\in \cG(\L_1,\L_2,\GG)\Big\}\nonumber\\
\cB(\GG_2,\GG)&:=&\Big\{(\Theta,\vartheta): \ |\Theta|< \frac{1}{2}\min\{\GG,\GG_2\}, |\vartheta|< \frac{\p}{2}\Big\}
\eeqa
with
\beqano
\cL(\GG)&:=&\Big\{\L=(\L_1,\L_2): \ \L\in \cL,\quad \L_1>\GG+\frac{2}{ c}\sqrt{\a_+}\L_2\Big\}\nonumber\\
\cG(\L_1,\L_2,\GG)&:=&\Big(\GG_-,\GG_+\Big),\qquad \GG_-:=\frac{2}{ c}\sqrt{\a_+} \L_2\qquad \GG_+:=\min\Big\{\L_1-\GG, \L_2\Big\}.\nonumber\\
\eeqano
where $\cL$ is as in~\equ{L0}, while $c$ is an arbitrarily fixed number in $(0,1)$.
We need to establish two kinds of conditions.

\paragraph{Geometric conditions}

First of all, we need that  the planets' eccentricities $e_1$, $e_2$ stay strictly confined in $(0,1)$. Then the following inequalities are to be satisfied: 
\beq{C1C2}0<\|{\rm C}^\ppu_{{\rm p}}\|<\L_1\qquad 0<\GG_2<\L_2,\eeq
with $\|{\rm C}^\ppu_{{\rm p}}\|$ as in~\equ{C}.
As remarked above, $\|{\rm C}^\ppu_{{\rm p}}\|$ may vanish only for
$$ (\GG_2,\vartheta)= (\GG,\p). $$

Since we deal with the equilibrium $(\uparrow\downarrow)$ (which holds for $(\Theta,\vartheta)=(0,0)$), the occurrence of this equality  is  automatically excluded, limiting the values of the coordinates $(\Theta,\vartheta)$ in the set $\cB$ in~\equ{B} since in this case
\beq{C1lowerbound}\|{\rm C}^\ppu_{{\rm p}}\|^2\ge \frac{3}{4}\GG^2.\eeq

Moreover, the two right inequalities in~\equ{C1C2} are satisfied taking 
\beqno\GG_2<\min\Big\{\L_1-\GG, \L_2\Big\}=\GG_+\eeqno
where we have used the triangular inequality $\|{\rm C}^\ppu_{{\rm p}}\|\le \|{\rm C}_{{\rm p}}\|+\|{\rm C}^\ppd_{{\rm p}}\|=\GG+\GG_2$.
\paragraph{Non-collision conditions} 

We have to exclude possible encounters of the planets with the sun and  each other.
Collisions of the inner planet with the sun are excluded by~\equ{B}. Indeed, using~\equ{C1lowerbound}, with $\L_1:=k_+\L_2$,
$$1-e_1^2=\frac{\|{\rm C}^\ppu_{{\rm p}}\|^2}{\L_1^2}\ge \frac{3}{4}\frac{\GG^2}{(\L_1)^2}$$
whence the minimum distance of the inner planet with the sun $a_1(1-e_1)$ is positive.
In order to avoid planetary collisions, it is typical to ensure the following inequality: 
$$a_1(1+e_1)<{c^2}a_2(1-e_2)$$
with $0<c<1$.
A sufficient condition for it is
	\beqno\GG_2\ge \frac{2}{{c}}\sqrt{\a_+} \L_2=\GG_-.\eeqno
 Indeed, if this inequality is satisfied, one has
$$a_1(1+e_1)<2a_1<\frac{a_2}{2}\frac{\GG_2^2{c^2}}{\L_2^2}=\frac{a_2}{2}(1-e_2^2){c^2}<a_2(1-e_2){c^2}.$$

\subsubsection{General properties of the secular problem}
{ We  call {\it Kepler map}  any canonical change of coordinates
\beq{Kepler maps}\textrm{\rm k}=(\L_1, \L_2,\ell_1, \ell_2,u,v)\in \cL\times\torus^2\times V\to (y^\ppu_{{\rm  k}}, y^\ppd_{{\rm  k}}, x^\ppu_{{\rm  k}}, x^\ppd_{{\rm  k}})\in \real^{12}\eeq
 such that 
 \beqno
\frac{|y^\ppj_{{\rm  k}}|^2}{2{\rm  m}_j}-\frac{{\rm  m}_j{\rm  M}_j}{|x_{{\rm  k}}^\ppj|}={\rm h}_{{\rm  k}}^\ppj(\L_j)\qquad j=1, 2,\eeqno
  where   $\cL\subset \real^2$, $V\subset\real^8$ are  open and connected sets, {and $\hh_\kk^\ppj$ is as in~\equ{h(L)}}. {Different Kepler maps are related by canonical changes}}
\vskip.1in
  
{\beqno\textrm{\rm k}=(\L_1, \L_2,\ell_1, \ell_2,u,v)\to \textrm{\rm k}'=(\L_1, \L_2,\ell'_1, \ell'_2,u',v')\eeqno
   which leave the $\L_i$'s unvaried. In terms of any Kepler map
the Hamiltonian~\equ{helio} takes the aspect
\beqno{\rm H}_{{\rm  k}}={\rm h}_{{\rm  k}}(\L_1,\L_2)+\m f_{{\rm  k}}(\L_1,\L_2,\ell_1,\ell_2,u,v)\eeqno  where}
  
{\beq{kepler}
{\rm h}_{{\rm  k}}(\L_1,\L_2)=-\frac{{\rm  m}_1^3{\rm  M}_1^2}{2\L_1^2}-\frac{{\rm  m}_2^3{\rm  M}_2^2}{2\L_2^2},\qquad f_{{\rm  k}}=-\frac{m_1m_2}{|x_{{\rm  k}}^\ppu-x_{{\rm  k}}^\ppd|}+\frac{y_{{\rm  k}}^\ppu\cdot y_{{\rm  k}}^\ppd}{m_0}.
\eeq}

If ${\rm k}$ is any  Kepler map, the  ``secular problem'' is the $(\ell_1, \ell_2)$--independent system
\beqa{secular system}
\HH^{\rm av}=\hh_{\rm k}+\m f^{\rm av}_{\rm k}
\eeqa
with
\beqa{av} f^{\rm av}_{\rm k}:=\frac{1}{(2\p)^2}\int_{[0,2\p]^2}f_{\rm k}d\ell_1d\ell_2\ .\eeqa
Some properties are enjoyed by all such systems, and are listed below.
\paragraph{The indirect part  in the perturbing function~\equ{kepler} does not contribute to  $f^{\rm av}_{\rm k}$.} Indeed, for any Kepler map, 
the impulses $y^\ppj_{\rm k}$ satisfy
$$y_{\rm k}^\ppj=\frac{\mm_j^2\MM_j}{\L_j}\partial_{\ell_j}\widehat x^\ppj_{\rm k}\ ,\quad {\rm with}\quad \widehat x^\ppj_{\rm k}:=\frac{x^\ppj_{\rm k}}{a_j}\ .$$

\paragraph{Expanding the averaged Newtonian potential in terms of the semi--axes ratio}
\beq{expansion}{f^{\rm av}_{\rm k}}=-\frac{m_1 m_2}{a_2}\Big({f^{\rm av}_{\rm k}}^\ppo+{f^{\rm av}_{\rm k}}^\ppu+\cdots\Big)\eeq
{\it then}
\beq{trivials}{f^{\rm av}_{\rm k}}^\ppo=1\qquad {f^{\rm av}_{\rm k}}^\ppu\equiv 0.\eeq
Indeed, ${f^{\rm av}_{\rm k}}^\ppo$ is the averaged Keplerian potential, given by $\frac{1}{a_2}$, while 
${f^{\rm av}_{\rm k}}^\ppu$ is the average of $x^\ppu_{\rm k}\cdot \frac{x^\ppd_{\rm k}}{\|x^\ppu_{\rm k}\|^3}$, which vanishes, because the second term is proportional to $\partial_{\ell_2}y^\ppd_{\rm k}$.
\paragraph{The term of order 2 in the expansion~\equ{expansion} admits $\GG_2$ as a first integral.} This circumstance has been pointed out in~Ref.\cite{harrington69}. In particular, in all Kepler maps such that $\GG_2$ is an action, ${f^{\rm av}_{\rm k}}^\ppd$ does not depend on its conjugate variable, and hence depends on one conjugated couple only. 
Therefore, the secular, truncated system
\beqa{averaged truncated}
\hh_{\rm eff}=\hh_{\rm k}
-\m\frac{\mm_1 m_2}{a_2}\Big(1+{f^{\rm av}_{\rm k}}^\ppd\Big)
\eeqa
which is obtained from~\equ{secular system} by truncating terms of order ${\rm O}(\m\a^3)$,
is integrable and, in particular, {\it one--dimensional}.
This property holds true, in particular, for ${\rm jrd}$ and ${\rm p}$. In the case of ${\rm jrd}$, ${f^{\rm av}_{\rm jrd}}^\ppd$, depends on the couple $(\GG_1,\g_1)$, whose bifurcation diagram has been studied, e.g.,  in~Refs.\cite{ferrerO94, lidovZ76, harrington69}. In the case of ${\rm p}$, ${f^{\rm av}_{\rm p}}^\ppd$, depends on the couple $(\Theta,\vartheta)$ and, due to  the equilibria~\equ{equilibria}, the dynamics of~\equ{averaged truncated} can be studied by {\it convergent} Birkhoff series.
More in general, the following formula has been proved in~Ref.\cite[Appendix B]{pinzari15}:
\beqa{fav2}
{f_\textrm{\rm k}
^{\rm av}}^\ppd&=&-\frac{\a^2}{8}\frac{\L_{2}^3}{\L_1^2\|{\rm C}_{{\rm  k}}^\ppd\|^5}
\Big[
5(3({\rm P}^\ppu_{{\rm  k}}\cdot {\rm C}^\ppd_{{\rm  k}})^2-\|{\rm C}^\ppd_{{\rm  k}}\|^2)\L_1^2\nonumber\\
&&\qquad-\; 3(4({\rm P}^\ppu_{{\rm  k}}\cdot {\rm C}^\ppd_{{\rm  k}})^2-\|{\rm C}^\ppd_{{\rm  k}}\|^2)\|{\rm C}^\ppu_{{\rm  k}}\|^2+3
({\rm C}_{{\rm  k}}^\ppu\times {\rm C}_{{\rm  k}}^\ppd\cdot {\rm P}_{{\rm  k}}^\ppu)^2
\Big].
\eeqa

\subsubsection{\label{bifurcations}Instability features in the secular problem}
In the case of ${\textrm{\rm k}}={\rm p}$, Equation~\equ{fav2} gives 
\beqa{f2}
{f_{\rm p}
^{\rm av}}^\ppd
= \a^2\big(\PP_0+\PP_1\big)
\eeqa
with
\beqa{f0P}
{\rm P}_0&:=&-\frac{1}{8}\frac{\L_{2}^3}{\L_1^2\GG_2^3}\Big(-5\L_1^2+3(\GG+\GG_2)^2\Big)\nonumber\\
{\rm P}_1&:=&-\frac{1}{8}\frac{\L_{2}^3}{\L_1^2\GG_2^5}\Big[
15\L_1^2\Theta^2-3(4\Theta^2-\GG_2^2)\Big(\GG^2+\GG_2^2-2\Theta^2+2\sqrt{(\GG_2^2-\Theta^2)(\GG^2-\Theta^2)}\cos{\vartheta}\Big)\nonumber\\
&&\qquad-\; 3\GG_2^2(\GG+\GG_2)^2+3(\GG_2^2-\Theta^2)(\GG^2-\Theta^2)\sin^2{\vartheta}\Big].
\eeqa
 Note that ${f_{\rm p}
^{\rm av}}^\ppd$ has been split so that $\a^2{\rm P}_0={f_{\rm p}
^{\rm av}}^\ppd\big|_{(\Theta,\vartheta)=(0,0)}$, whence $\PP_1$ vanishes for $(\Theta,\vartheta)=(0,0)$. In the next proposition, we study the properties of $\PP_1$.

Recall the definition of the sets $\cA(\GG)$, $\cB(\GG_2,\GG)$ in~\equ{B}.

\begin{proposition}\label{maybehyperbolic}
For any fixed $\GG\in \real_+$, there exist $\cA_{\rm u}(\GG)\subset\cA(\GG)$ and a neighborhood $\cB_{\rm u}$ of  $(0,0)$,
with $\cB_{\rm u}\subset \cB(\GG_2,\GG)$ for all $\GG_2\in {\cal G}$,
such that $\HH_{\rm p}$ is real--analytic on
\beqa{DU}\cD_{\rm u}:=\cA_{\rm u}\times \cB_{\rm u}\times\torus^3\eeqa
and $(\Theta,\vartheta)=(0,0)$ is a hyperbolic  equilibrium point for $\PP_1$. More precisely, there exist two functions  $\o$, $\O$ of $\L_1$, $\L_2$, $\GG_2$ and the parameter $\GG$, 
with $\o>0$ such that, if
\beq{p0q0}p_0:=\frac{\Theta-\o \vartheta}{\sqrt{2\o}},\qquad q_0:=\frac{\Theta+\o \vartheta}{\sqrt{2\o}}\eeq
one has
\beqno{\rm P}_1=\O p_0 q_0+{\rm O}(p_0,q_0; \L_2,\L_2,\GG_2,\GG)^4.\eeqno
\end{proposition}

We shall prove the proposition with
\beqa{assumptions0}
\O&:=&-\frac{3}{4}\frac{\L_{2}^3}{\L_1^2\GG_2^4}\sqrt{\big(5\L_1^2\GG -(\GG+\GG_2)^2 (4 \GG+\GG_2)\big)\big(\GG_2- \GG\big)}\nonumber\\
\o&:=&\GG\GG_2\sqrt{\frac{\GG_2- \GG}{5\L_1^2\GG -(\GG+\GG_2)^2 (4 \GG+\GG_2)}}\nonumber\\
\cA_{\rm u}(\GG)&:=&\Big\{(\L_1,\L_2)\in {\cal L}_{\rm u}(\GG),\quad \GG_2\in{\cal G}_{\rm u}(\L_1,\L_2,\GG)\Big\}\nonumber\\
\cB_{\rm u}(\GG)&:=&\Big\{(\Theta,\vartheta): \ |\Theta|< \frac{\GG}{2}, |\vartheta|< \frac{\p}{2}\Big\}
\eeqa 
where
\beqa{assumptions}\cL_{\rm u}(\GG)&:=&\Big\{\L=(\L_1,\L_2)\in \cL: \ 5\L_1^2\GG -(\GG+\frac{2}{ c}\sqrt{\a_+}\L_1)^2 (4 \GG+\frac{2}{ c}\sqrt{\a_+}\L_1)>0,\nonumber\\
&& \hspace*{10em} \L_2>\GG, \L_1>\max\{\GG+\frac{2}{ c}\sqrt{\a_+}\L_2, 2\GG\}\Big\}\nonumber\\
\cG_{\rm u}(\L_1,\L_2,\GG)&:=&\Big(\ovl\GG_-, \ovl\GG_+\Big)
\eeqa
where $\cL$ is as in~\equ{L0} and, if $\GG^{\star}(\L_1,\GG)$ is the unique positive root of the cubic polynomial $\GG_2\to 5\L_1^2\GG -(\GG+\GG_2)^2 (4 \GG+\GG_2)$, then

\beq{assumptions3}\ovl\GG_-:=\max\{\frac{2}{ c}\sqrt{\a_+}\L_2, \GG\}\qquad \ovl\GG_+:=\min\{ \L_2, \GG^{\star}\}.\eeq
Implicitly, we shall prove that
\beq{assumptions7}\ovl\GG_-< \ovl\GG_+.\eeq

{\bf Proof.}
The expansion of $\PP_{1}$ in~\equ{f0P} around $(0,0)$ is 
 \beqano
\PP_{1} = -\frac{1}{8}\frac{\L_{2}^3}{\L_1^2\GG_2^5} \times \Big[\frac{3}{\GG}a(\L_1,\GG_2;\GG)\Theta^2+3\GG\GG_2^2b(\GG_2;\GG)
\vartheta^2+o_2(\Theta,\vartheta)\Big]
\eeqano
where

\beq{a b}a(\L_1,\GG_2;\GG):=5\L_1^2\GG -(\GG+\GG_2)^2 (4 \GG+\GG_2)\quad {\rm and}\quad b(\GG_2;\GG):=\GG- \GG_2.\eeq
Both $\GG_2\to a(\L_1,\GG_2;\GG)$ and $\GG_2\to b(\GG_2;\GG)$, as functions of $\GG_2$ decrease monotonically from a positive value (respectively, $\GG(5\L_1^2-4\GG^2)$ and $\GG$) to $-\infty$ as $\GG_2$ increases from $\GG_2=0$ to $\GG_2=+\infty$. The function $a(\L_1,\GG_2;\GG)$ changes its sign for $\GG_2$ equal to a suitable unique positive value $\GG^{\star}(\L_1,\GG)$, while  $b(\GG_2;\GG)$ does it for $\GG_2=\GG$. We note that (i) inequality $\GG<\min\{\GG_+, \GG^{\star}\}$ follows immediately from the assumptions~\equ{assumptions} (in particular, the two last ones) and (ii), more generally, that
 $\GG^{\star}\leq\GG$ is equivalent to $\L_1\leq2\GG$. Since, for our purposes, we have to exclude $\GG^{\star}=\GG$ (otherwise, $a(\L_1,\GG_2;\GG)$ and $b(\GG_2;\GG)$ would be simultaneously positive and simultaneously negative, and no hyperbolicity would be possible), we distinguish two cases. 
 
 \begin{itemize}
 \item[(a)] $\GG>\frac{2}{ c}\sqrt{\a_+} \L_2$ and $\GG+\frac{2}{ c}\sqrt{\a_+} \L_2<\L_1<2\GG$. In this case $\GG^\star<\GG$. We show that no such $\cG_{\rm u}$ can exist in this case. In fact, since $\GG^{\star}<\GG$, in order that the interval $(\GG^{\star}, \GG)$ and the set  ${\cal G}$ have a non-empty intersection,  one should have, necessarily, $\GG_+=\sup{\cal G}>\GG^{\star}$, hence, in particular,  $\L_1-\GG>\GG^{\star}$. Using the definition of $\GG^{\star}$, this would imply $\L_1>2\GG$, which is a contradiction.

 \item[(b)] $\L_1>\max\{2\GG, \GG+\frac{2}{ c}\sqrt{\a_+} \L_2\}$. In this case $\GG<\GG^{\star}<\L_1-\GG$. In order that the interval $(\GG, \GG^{\star})$ and the set  ${\cal G}$ have a non-empty intersection, we need
\beq{conditions}\GG_-<\GG^{\star}\qquad {\rm and}\qquad \GG_+>\GG\eeq and such intersection will be given by the interval $\cG_{\rm u}$ as in~\equ{assumptions}. Note that the definition of $\ovl\GG_+$ does not include $\L_1-\GG$ in the brackets because, as noted, $\GG^{\star}<\L_1-\GG$. 
But~\equ{conditions} are equivalent to
\equ{assumptions}. 
\end{itemize}
\begin{remark}\label{elliptic}\rm
The ``bifurcation'' towards the hyperbolic behavior that Proposition~\ref{maybehyperbolic} talks about does not appear in the case of the equilibria $(\uparrow\uparrow)$ and $(\downarrow\uparrow)$, which, in contrast, are always
  {\it  elliptic}. Indeed, in such cases, one obtains an expansion analogous to~\equ{f0P}, with the coefficients $a$, $b$ in~\equ{a b}  to be replaced by
$$
\widehat a=5\L_1^2\GG -(\GG-\GG_2)^2 (4 \GG-\GG_2),\qquad \widehat b= \GG+ \GG_2.
$$
Clearly,
 $\widehat b$ is positive for all $\GG_2$ and $\widehat a$ is so for $\GG_2\ge 4\GG$. On the other hand, when $\GG_2< 4\GG$, inequality $|\GG-\GG_2|<\L_1$  implies
 $$\widehat a\ge \L_1^2(\GG+\GG_2)>0.$$
 Therefore,  $\widehat a$ and  $\widehat b$ have always the same (positive) sign. This circumstance has been worked out in~Ref.\cite{pinzari15}, in the more general situation with $N\ge 2$ planets, 
in order to infer the existence of quasi--periodic motions with maximal number of frequencies away from the constraint, we have talked about of in Sections~\ref{The canonical setting}--\ref{domain of rps+}, of small eccentricities and inclinations.
\end{remark}

\subsubsection{\label{conjectures}
Unperturbed hyperbolic motions
}
Let us consider the Hamiltonian $\HH_{\rm p}$ in~\equ{p Ham}, on the domain $\cD_{\rm u}$ in~\equ{DU}. 
Standard averaging theory (e.g.,~Ref.\cite[Normal Form Lemma]{poschel93}) allows to eliminate the dependence on $\ell_1$, $\ell_2$ at a higher order,  conjugating, via a real--analytic, $\m$--close to the identity, canonical transformation,  $\HH_{\rm p}$ to a new Hamiltonian, that we denote as
$$\HH'_{\rm p}=\hh_{\rm k}(\L_1,\L_2)+\m f_{\rm p}^{\rm av}(\L_1,\L_2,\GG_2,\Theta,{\rm g}_2,\vartheta;\GG)+\m^\s f'_{\rm p}(\L_1,\L_2,\GG_2,\Theta,\ell_1,\ell_2,{\rm g}_2,\vartheta;\GG)$$
where $f_{\rm p}^{\rm av}$ is as in~\equ{av}, while  $1<\s<2$. In view of
\equ{expansion}, we  can thus split $\HH'_{\rm p}$ as 
$$\HH'_{\rm p}=\hh_{\rm eff}+ f_{\rm eff}$$
where $\hh_{\rm eff}$ is as in~\equ{averaged truncated}, with k=p, namely,
\beqa{H0}\hh_{\rm eff}:= {\rm h}_{{\rm  k}}(\L_1,\L_2)-\m\frac{m_1 m_2}{a_2}\Big(1+\a^2({\rm P}_0+\PP_{1})\Big)\eeqa
while
 $f_{\rm eff}:=\m\a^3\widehat\PP+\m^\s f'_{\rm p}$ is a smaller remainder. The leading term,
$\hh_{\rm eff}$, possesses, by Proposition~\ref{maybehyperbolic}, a hyperbolic fixed point at $(\Theta,\vartheta)=(0,0)$
and hence a family of three-dimensional tori
with linear motions of $(\ell_1, \ell_2, {\rm g}_2)$,
 having equation
$$\cT_{\L_1^\star,\L_2^\star,\GG_2^\star}=\Big\{(\L_1,\L_2,\GG_2)=(\L_1^\star,\L_2^\star,\GG_2^\star),\quad (\ell_1,\ell_2,\textrm{\rm g}_2)\in \torus^3,\quad (\Theta,\vartheta)=(0,0)\Big\}$$
parametrized by $(\L_1^\star,\L_2^\star,\GG_2^\star)\in \cA_{\rm u}(\GG)$. Together to such tori, the unperturbed system also possesses
  two four-dimensional manifolds, the {\it whiskers}, \beqa{unperturbed whiskers}
&&{\cal W}_{\L_1^\star,\L_2^\star,\GG_2^\star,\varepsilon}^{\rm s,loc}=\Big\{(\L_1,\L_2,\GG_2)=(\L_1^\star,\L_2^\star,\GG_2^\star), (\ell_1,\ell_2,\textrm{\rm g}_2)\in \torus^3,\quad q_0=0, |p_0|<\varepsilon\Big\}\nonumber\\
&&{\cal W}_{\L_1^\star,\L_2^\star,\GG_2^\star,\varepsilon}^{\rm u,loc}=\Big\{(\L_1,\L_2,\GG_2)=(\L_1^\star,\L_2^\star,\GG_2^\star), (\ell_1,\ell_2,\textrm{\rm g}_2)\in \torus^3,\quad p_0=0, |q_0|<\varepsilon\Big\}
\eeqa
(with $p_0$, $q_0$ as in~\equ{p0q0})  including the tori $\cT_{\L_1^\star,\L_2^\star,\GG_2^\star}$  as subsets, with motions asymptotic to $\cT_{\L_1^\star,\L_2^\star,\GG_2^\star}$ in the future/past, respectively.
 It is a reasonable expectation (motivated by the well established {\sc kam} theory) that such structure (tori and related whiskers) is preserved in the whole system, at least for those $\cT_{\L_1^\star,\L_2^\star,\GG_2^\star}$'s whose frequencies $\o^\star:=\partial_{(\L^*_2,\L^*_2,\GG^*_2)}\hh_{\rm eff}\big|_{(\Theta,\vartheta)=(0,0)}$ are highly irrational (Diophantine). \section{\label{appendix D}{On the canonical character of Jacobi--Radau--Deprit and {\rm p}--coordinates}}

In this section we discuss, in a unified way, the canonical character of the {\rm jrd} coordinates in~\equ{J}--\equ{coordinates} and  the {\rm p}--coordinates in~\equ{belle*}.  In both cases, we reduce to the discussion to the case $N=2$, as needed in the paper. Both such sets of coordinates are defined for a general number of particles. The generalization of {\rm jrd} can be deduced from the original Deprit's paper,~Ref.\cite{deprit83}. See also ~Ref.\cite{pinzari-th09}, for a direct,  inductive approach. As for the generalization of {\rm p}, one can look at~Ref.\cite{pinzari15}.

The proof of the canonical character of {\rm jrd} and {\rm p} will be based on the Delaunay coordinates, that here we recall, and a simple lemma.
\begin{itemize}
\item[$\star$]   Delaunay coordinates (see~Ref.\cite{GG83}),  six for every body, here denoted as \beqa{d}\textrm{\rm d}=(\L_j,  \GG_j,\HH_j, \ell_j, \ovl{\rm g}_j,\hh_j)\qquad j=1,\ 2\eeqa
are defined as follows. The coordinates $\L_j$, $ \GG_j$, $\ell_j$ are as in~\equ{coordinates}, while, if \beq{Del nodes} n_j:=k^\ppt\times\CC^\ppj\eeq
 and $\a_w(u,v)$ as said  in Section~\ref{The canonical setting}, then 
\beq{Delaunay} \ovl{\rm g}_j=\a_{\CC^\ppj}(n_j,{\rm P}^\ppj)\ ,\qquad \hh_j=\a_{k^\ppt}(k^\ppu, n_j)\ ,\quad \HH_j=\CC^\ppj\cdot k^\ppt=\CC^\ppj_3\ ,\eeq
with $(k^\ppu, k^\ppd,k^\ppt)$ a prefixed orthonormal frame in $\real^3$.
\item[$\star$] {Let $R_e(\a)\in {\rm SO}(3)$ be a rotation by an angle $\a$ around the unit vector $e$; let
$$x(\a,\ovl x):=R_e(\a) \ovl x\ \forall\ \ovl x\in \real^3\ ,\quad \CC(x,y):=x\times y\ \forall\ x,\ y\in \real^3\ .$$
Then we have}
 \begin{lemma}\label{trivial}
{Given $\ovl y\in \real^3$ and varying $\ovl x$ and $\a$,}
$$ y\cdot d x=(\CC(\ovl x,\ovl y)\cdot e)d\a+\ovl y\cdot d\ovl x=(\CC(x,y)\cdot e)d\a+\ovl y\cdot d\ovl x
\ 
.$$
 \end{lemma}
 {{\bf Proof.} Varying $\a$ and letting $\ovl x$, $\ovl y$ fixed, 
\beqano
y(\a,\ovl y)\cdot d\big(R_e(\a)\ovl x\big)&=&\big(R_e(\a)\ovl y\big)\cdot d\big(R_e(\a)\ovl x\big)=\big(R_e(\a)\ovl y\big)\cdot\big(e\times R_e(\a)\ovl x\big)\nonumber\\
&=&\big(R_e(\a)\ovl x\big)\times \big(R_e(\a)\ovl y\big)\cdot e\nonumber\\
&=&C(\ovl x,\ovl y)\cdot e=C(x, y)\cdot e
\ . \eeqano
Letting also $\ovl x$ vary, one has the thesis.
 }
\end{itemize}
Now we proceed with proving the canonical character of {\rm jrd} and {\rm p}, {by showing that they are canonically related to the coordinates ${\rm d}$ in~\equ{d}.}
Since the couples $(\L_j, \ell_j)$ are in common to {\rm d}, {\rm jrd} and {\rm p}, {namely,
	\beqano
	{\rm d}=\big(\L_1,\L_2,\ell_1,\ell_2,\widehat{\rm d}\big)\ ,\quad {\rm jrd}=\big(\L_1,\L_2,\ell_1,\ell_2,\widehat{\rm jrd}\big)\ ,\quad {\rm p}=\big(\L_1,\L_2,\ell_1,\ell_2,\widehat{\rm p}\big)
	\eeqano
	with
	\beqano
	&&\widehat{\rm d}=(\HH_1,\HH_2,\GG_1,\GG_2,\hh_1,\hh_2,\ovl{\rm g}_1,\ovl{\rm g}_2)\ ,\quad \widehat{\rm jrd}=(\ZZ,\GG,\GG_1,\GG_2,\zeta,\g,\g_1,\g_2)\nonumber\\
	&& \widehat{\rm p}=(\ZZ,\GG,\Theta,\GG_2,\zeta,{\rm g},\vartheta,{\rm g}_2)
	\eeqano
}
{and the changes
\beq{2changes}\widehat{\rm d}\to \widehat{\rm jrd}\to \widehat{\rm p}\eeq
do not depend on $(\L_j, \ell_j)$,}
 we just need to
check that 
such changes {in~\equ{2changes}} are  canonical. We shall prove that 
\begin{theorem}\label{theo: homo} The changes of coordinates in~\equ{2changes} preserve the standard $1$--form:
\beq{homo}\sum_{i=1}^2(\HH_id\hh_i+\GG_id\ovl{\rm g}_i)=\ZZ d\zeta+\GG d\g+\GG_1d\g_1+\GG_2d\g_2=\ZZ d\zeta+\GG d{\rm g}+\Theta d\vartheta+\GG_2d{\rm g}_2\ .\eeq
\end{theorem}
We shall use many times the following definitions.
\begin{definition}
If $n\perp n'\in \real^3$, we  denote as $\FF\sim(n,\cdot,n')$ the orthonormal frame  $\FF=(\frac{n}{|n|}, \frac{n'\times n}{|n'\times n|},\frac{n'}{|n'|})$.
\end{definition}
\begin{definition}\label{good changes}
We denote as
$$\FF\to\FF'$$
any couple $(\FF,\FF')$ of orthonormal frames, with $\FF=(i,j,k)$, $\FF'=(i',j',k')$, such that $i'\parallel \pm k\times k'$.
\end{definition}
 Observe that, in such case, the transformation of coordinates
which relates the coordinates ${\rm X}'$ relatively to $\FF'$ to  the coordinates ${\rm X}$ relatively to $\FF$  is
\beq{general case}{\rm X}=R_3(\psi)R_1(\iota){\rm X}'\eeq
where {$R_1$, $R_3$ are as in~\equ{R1R3};} $\iota$, the ``mutual inclination between $\FF$ and $\FF'$'',  is the convex angle  between $k$ and $k'$, while $\psi$, called ``longitude of the node of $\FF'$ with respect to $\FF$, is defined by $\psi:=\a_k(i,i')$.

{\bf Proof.}
Let  $\FF_0$ be a prefixed reference frame, and let $\FF^\ppj_{\rm d}\sim(n_j,\cdot, \CC^\ppj)$, with $n_j$ as in~\equ{Del nodes}. Then we have 
\beq{Del rotations}\FF_0\to\FF_{\rm d}^\ppj\ .\eeq
Let  $\PP^\ppj$ be the coordinates of the $j^{\rm th}$ perihelion with respect to  $\FF_0$, and denote as ${\rm Q}^\ppj:=\widehat\CC^\ppj\times\PP^\ppj$, with 
 $\widehat\CC^\ppj=\frac{\CC^\ppj}{|\CC^\ppj|}$, so  that $(\PP^\ppj, {\rm Q}^\ppj,\widehat\CC^\ppj)$ is an orthonormal triple in $\real^3$.
The coordinates of such vectors  relatively to $\FF^\ppj_{\rm d}$,  are
$$\PP^\ppj_{\rm d}=R_3(\ovl{\rm g}_j)e_1\qquad {\rm Q}^\ppj_{\rm d}=R_3(\ovl{\rm g}_j)e_2\ .$$
Therefore, by~\equ{general case},~\equ{Del rotations} and the definitions in~\equ{Delaunay}, we have
 \beqno\PP^\ppj=R_3(\hh_j)R_1(\ovl i_j)R_3(\ovl{\rm g}_j)e_1\ ,\qquad {\rm Q}^\ppj=R_3(\hh_j)R_1(\ovl i_j)R_3(\ovl{\rm g}_j)e_2\ ,\eeqno
 where $\cos\ovl i_j=\frac{\HH_j}{\GG_j}$.
 Then in view of  Lemma~\ref{trivial}, we obtain (using $e_1\cdot R_3(\g_j)(e_1\times e_2)=e_1\cdot e_3=0$ and $\PP^\ppj\times {\rm Q}^\ppj=\widehat\CC^\ppj$)
  $${\rm Q}^\ppj\cdot d\PP^\ppj=\widehat\CC^\ppj\cdot e_3 d\hh_j+d\ovl{\rm g}_j\ .$$
   Multiplying by $\GG_j=|\CC^\ppj|$ and recognizing that 
  $\GG_j\widehat\CC^\ppj\cdot e_3=\HH_j$, we then have
  \beq{DelPQ}\sum_{i=1}^2(\HH_id\hh_i+\GG_id\ovl{\rm g}_i)=\sum_{j=1}^2|\CC^\ppj|{\rm Q}^\ppj\cdot d\PP^\ppj\ .\eeq
  Now we compute the right hand side of this equation, using the {\rm jrd} and {\rm p}--coordinates. To this end, we need to express $\PP^\ppj$ and ${\rm Q}^\ppj$ in terms of such two sets. To accomplish this, we observe that, in the sense of Definition~\ref{good changes}, 
  
  \begin{itemize}
  \item[$\star$] In the case of {\rm jrd}, we have the  ``tree'' of changes of frames, 
  \beqa{tree}\begin{array}
    {llllllllllllllllllll}{\rm F}_0&\to&\FF_*\to{\rm F}_{{\rm  jrd}}^\ppu\\
    &&\downarrow\\
    &&{\rm F}_{{\rm  jrd}}^\ppd
\end{array}
\eeqa
where  ${\rm F}_{{\rm  jrd}}^\ppj\sim(\n,\cdot,\CC^\ppj)$, while $\FF_*\sim(\n_1,\cdot,\CC)$ is the {\it invariable frame}.
  \item[$\star$] In the case of {\rm p}, we have the ``chain'' 
\beq{chain}\FF_0\to \FF_*\to \GG_{{\rm p}}^\ppu\to \FF_{{\rm p}}^\ppd\to \GG_{{\rm p}}^\ppd\eeq
 where $\FF_0$, $\FF_*$ are as in the previous item, while
  $$\GG_{{\rm p}}^\ppu\sim({\rm n}_1,\cdot, \PP^\ppu)\ ,\qquad \FF^\ppd_{{\rm p}}\sim (\n_2,\cdot,\CC^\ppd)\ ,\qquad \GG_{{\rm p}}^\ppd\sim({\rm n}_2,\cdot,\PP^\ppd)\ .$$
   \end{itemize}
   Therefore,
\begin{itemize}
\item[$\star$] Recognizing (by the analysis of the triangle formed by $\CC^\ppu$, $\CC^\ppd$ and $\CC=\CC^\ppu+\CC^\ppd$)
 that the  inclinations $i$, $i_1$, $i_2$ between
  $\FF_0$ and $\FF_*$,$\FF_*$ and $\FF^\ppu_{{\rm  jrd}}$, $\FF_*$ and $\FF^\ppd_{{\rm  jrd}}$ are given by~\equ{D incli},  while the longitudes of the nodes  are, respectively, $\zeta$,  $\g$, $\g+\p$, we find the formulae
  \beqa{generic}
  &&\PP^\ppj=R_3(\zeta)R_1(i)R_3(\g)R_1(s_ji_j)R_3(\g_j)e_1\nonumber\\
  && {\rm Q}^\ppj=R_3(\zeta)R_1(i)R_3(\g)R_1(s_ji_j)R_3(\g_j)e_2\eeqa
  where $s_1=-s_2=1$. {Here we have used $R_3(\p)R_1(\a)=R_1(-\a)R_3(\p)$ and the definitions of the angles $\g_j$ in~\equ{coordinates}. Observe, incidentally, that
  the formulae in~\equ{generic}  allow us to obtain the formulae in~\equ{yx}, since, as known, $y^\ppj$, $x^\ppj$ are related to   
$a_j$, ${\rm e}_j$,   $\PP^\ppj$ and ${\rm Q}^\ppj$ via the classical relations in~\equ{xyP}.
with $a_j$ as in~\equ{aL}; ${\rm e}_j$, $\zeta_j$ as in~\equ{ej},~\equ{zetaj}.
  }
\item[$\star$] The expressions of $\PP^\ppj$, ${\rm Q}^\ppj$ in terms of the ${\rm p}$ coordinates have been given in~\equ{P formulae},~\equ{Q1Q2}.
\end{itemize}
We are now ready to compute the right hand side of~\equ{DelPQ}, in terms of {\rm jrd} and {\rm p}. To this scope, we shall use Lemma~\ref{trivial}.
  \begin{itemize}
  \item[$\star$] 
Using the formulae in~\equ{generic}, iterate applications of Lemma~\ref{trivial} and linear algebra, we obtain
  $${\rm Q}^\ppj\cdot d\PP^\ppj=\widehat\CC^\ppj\cdot e_3d\zeta+\widehat\CC^\ppj\cdot \n d i+\widehat\CC^\ppj\cdot kd\g+{\rm f}_j \cdot e_1di_j+ d\g_j$$
with $\n:=R_3(\zeta)e_1=(\cos\zeta,\sin\zeta, 0)$, $k:=R_3(\zeta)R_1(\iota)e_3$, ${\rm f}_j:=R_1(s_j i_j)e_3$. Multiplying by $\GG_j=|\CC^\ppj|$, summing over $j=1$, $2$ and recognizing that $k$ has the direction of $\CC$,  $\n$ is orthogonal to $\CC$ and $\GG_1{\rm f}_1+\GG_2{\rm f}_2=\GG e_3$, we immediately obtain, after some cancellation,
\beq{homo2}\sum_{j=1}^\ppd|\CC^\ppj|{\rm Q}^\ppj\cdot d\PP^\ppj=\ZZ d\zeta+\GG d\g+\GG_1d\g_1+\GG_2d\g_2\ .\eeq  \item[$\star$] Using the formulae in~\equ{P formulae}, defining  $\widehat\CC^\ppu_1$,  $\widehat\CC^\ppu_2$,  $\widehat\CC^\ppd_1$,  $\widehat\CC^\ppd_2$, $\widehat\CC^\ppd_3$ via
\beqano
\widehat\CC^\ppu&=&R_3(\zeta)R_1(i)\widehat\CC^\ppu_1=R_3(\zeta)R_1(i)R_3({\rm g})R_1(i_1)\widehat\CC^\ppu_2\nonumber\\
\widehat\CC^\ppd&=&R_3(\zeta)R_1(i)\widehat\CC^\ppd_1=R_3(\zeta)R_1(i)R_3({\rm g})R_1(i_1)\widehat\CC^\ppd_2\nonumber\\
&=&R_3(\zeta)R_1(i) R_3({\rm g})R_1(i_1)R_3(\vartheta)R_1(i_2)\widehat\CC^\ppd_3
\eeqano
and applying iteratively  Lemma~\ref{trivial}, we obtain
\beqano
{\rm Q}^\ppu\cdot d\PP^\ppu&=&\widehat\CC^\ppu\cdot e_3d\zeta+\widehat\CC^\ppu_1\cdot e_1 di+\widehat\CC^\ppu_1\cdot e_3 d{\rm g}+\widehat\CC^\ppu_2\cdot e_1 di_1\nonumber\\
{\rm Q}^\ppd\cdot d\PP^\ppd&=&\widehat\CC^\ppd\cdot e_3d\zeta+\widehat\CC^\ppd_1\cdot e_1 di+\widehat\CC^\ppd_1\cdot e_3 d{\rm g}+\widehat\CC^\ppd_2\cdot e_1 di_1\nonumber\\
&+&\widehat\CC^\ppd_2\cdot e_3 d\vartheta+\widehat\CC^\ppd_3\cdot e_1 di_2+\widehat\CC^\ppd_3\cdot e_3 d{\rm g}_2
\eeqano
We multiply, as above, the first equation by $|\CC^\ppu|$, the second by $|\CC^\ppd|$, and take the sum of the two.
The sum of the first three respective terms gives, analogously to the previous item, $\ZZ d\zeta+\GG d{\rm g}$.
As for the remaining terms, we recognize that $\widehat\CC^\ppd_3=e_3$ so that $\widehat\CC^\ppd_3\cdot e_1=0$, $\widehat\CC^\ppd_3\cdot e_3=1$; $\widehat\CC^\ppd_2\cdot e_3=\Theta$; $|\CC^\ppu|\widehat\CC^\ppu_2\cdot e_1+|\CC^\ppd|\widehat\CC^\ppd_2\cdot e_1=(\CC^\ppu+\CC^\ppd)\cdot R_3(\zeta)R_1(i)R_3({\rm g})R_1(i_1) e_1=\CC\cdot R_3(\zeta)R_1(i)R_3({\rm g})R_1(i_1) e_1=\GG e_3\cdot R_3({\rm g})e_1=0
$. We finally obtain
\beq{homo3}\sum_{j=1}^2|\CC^\ppj|{\rm Q}^\ppj\cdot d\PP^\ppj=\ZZ d\zeta+\GG d{\rm g}+\Theta d\vartheta+\GG_2d{\rm g}_2\ .\eeq
  \end{itemize}
The collection of~\equ{DelPQ},~\equ{homo2} and~\equ{homo3} proves Theorem~\ref{theo: homo}.


\appendix

\section{}
 
\subsection{\label{lemma sui segni}Proof of Proposition~\ref{signs1} and other technicalities}
The starting point is the analytical expression of $\phi_{\rm jrd}$, as presented in Refs.\cite{pinzari-th09, chierchiaPi11b}, that here we recall.
\subsubsection{Analytical expression of $\phi_{\rm jrd}$}\label{expression of jrd}
{Let  ${\rm e}_j\in (0,1)$ be the eccentricity of the Keplerian orbit, as known, related to $\L_j$, $\GG_j$ via
\beq{ej}{\rm e}_j=\sqrt{1-\frac{\GG_j^2}{\L_j^2}}\eeq
$\zeta_j({\rm e}_j, \ell_j)$ the {\it eccentric anomaly}, defined 
as the unique solution of {\it Kepler equation}
\beq{zetaj}
\zeta_j-{\rm e}_j\sin\zeta_j=\ell_j\qquad j=1, \ 2
\eeq
and let $i$, $i_1$, $i_2$ be defined via
\beqa{D incli} i=\cos^{-1}\frac{\ZZ}{\GG}\ ,\quad  i_1=\cos^{-1}\frac{\GG_1^2+\GG^2-\GG_{2}^2}{2\GG\GG_1}\ ,\quad   i_{2}=\cos^{-1}\frac{\GG_2^2+\GG^2-\GG_1^2}{2\GG\GG_2}\ .\eeqa
By~\equ{Dep conditions}, $i$, $i_1$ and $i_2$ take values on $(0,\p)$.
Let, finally, $R_1$, $R_3$ denote the matrices
\beqa{R1R3}
&&R_1(\a):=\left(
\begin{array}{ccc}
1&0&0\\
0&\cos \a&-\sin\a\\
0&\sin\a&\cos\a
\end{array}
\right)\quad  R_3(\a):=\left(
\begin{array}{ccc}
\cos\a&-\sin\a&0\\
\sin\a&\cos\a&0\\
0&0&1
\end{array}
\right)\ .\eeqa
The map $\phi_{\rm jrd}$ in~\equ{jrd map} is defined as
\beqa{yx}
\left\{
\begin{array}{ll}
x_{\rm jrd}^\ppj=R_3(\zeta)R_1(i)R_3(\g)R_1(s_ji_j)R_3(\g_j)x^\ppj_{\rm orb}(\L_j,\GG_j,\ell_j)\\
y_{\rm jrd}^\ppj= R_3(\zeta)R_1(i)R_3(\g)R_1(s_ji_j)R_3(\g_j)y^\ppj_{\rm orb}(\L_j,\GG_j,\ell_j)
\end{array}
\right.
\eeqa
with $s_1=-s_2=1$ and
\beqa{xy orb}
\left\{
\begin{array}{ll}
x_{\rm orb}^\ppj=\frac{\L_j^2}{\mm_j^2\MM_j} \left(
\begin{array}{ccc}
\cos\zeta_j-{\rm e}_j\\
\sqrt{1-{\rm e}_j^2}\sin\zeta_j\\
0
\end{array}
\right)\\
y_{\rm orb}^\ppj=\frac{\mm_j^2\MM_j}{\L_j}\frac{1}{1-{\rm e}_j\cos\zeta_j} \left(
\begin{array}{ccc}
-\sin\zeta_j\\
\sqrt{1-{\rm e}_j^2}\cos\zeta_j\\
0
\end{array}
\right)
\end{array}
\right.
\eeqa
The~\equ{yx}
simply describe the successive rotations necessary to transform the coordinates of $x^\ppj$, $y^\ppj$ relatively to the ``orbital frame''  in \equ{xy orb}  into the coordinates relatively to the prefixed frame. Roughly, the choice $s_1=-s_2=1$ reflects the classical ``opposition of the nodes'' in the frame where the total angular momentum is vertical, i.e., the identity $\CC\times \CC^\ppu+\CC\times\CC^\ppd=\CC\times\CC\equiv0$. }

\subsubsection{Analytical expression  of $\phi^\complex_{\rm rps}$}
The formulae  of the $\phi^\complex_{\rm rps}$ map  may  be recovered rewriting~\equ{yx} in the form
\beqa{new formulae*}
\left\{
\begin{array}{ll}
x_{\rm jrd}^\ppj=\widetilde\RR_0\widetilde\RR_j\widetilde{\rm x}^\ppj_{\rm pl}\\
y_{\rm jrd}^\ppj=\widetilde\RR_0\widetilde\RR_j\widetilde{\rm y}^\ppj_{\rm pl}
\end{array}
\right.
\eeqa
where
\beqano
&&\widetilde\RR_0:=\widetilde R_{313}(\zeta,i)\ ,\quad \widetilde\RR_j:=\widetilde R_{313}(\zeta+\g, s_ji_j)\ ,\quad \widetilde{\rm x}^\ppj_{\rm pl}:=R_3(\zeta+\g+\g_j)x^\ppj_{\rm orb}(\L_j,\GG_j,\ell_j)\nonumber\\
&& \widetilde{\rm y}^\ppj_{\rm pl}:=R_3(\zeta+\g+\g_j)y^\ppj_{\rm orb}(\L_j,\GG_j,\ell_j)
\eeqano
with $ x^\ppj_{\rm orb}$,$ y^\ppj_{\rm orb}$ as in~\equ{xy orb},
\beqa{formulae}
\widetilde R_{313}(\a,\b)&:=&R_3(\a)R_1(\b)R_3(-\a)\nonumber\\
&=&\left(
\begin{array}{ccc}
1-\sin^2\a(1-\cos{\b})&\sin{\a}\cos{\a}(1-\cos{\b})&\sin{\a}\sin{\b}\\
\sin{\a}\cos{\a}(1-\cos{\b})&1-\cos^2\a(1-\cos{\b})&-\cos{\a}\sin{\b}\\
-\sin{\a}\sin{\b}&\cos{\a}\sin{\b}&\cos{\b}
\end{array}
\right)
\eeqa
Note that,
{in particular,     the matrices $\widetilde\RR_0$, $\widetilde\RR_1$, $\widetilde\RR_2$ 
    are given by
    {\small    $$\widetilde\RR_0=\left(
\begin{array}{ccc}
1-\sin^2\zeta(1-\cos{i})&\sin{\zeta}\cos{\zeta}(1-\cos{i})&\sin{\zeta}\sin{i}\\
\sin{\zeta}\cos{\zeta}(1-\cos{i})&1-\cos^2\zeta(1-\cos{i})&-\cos{\zeta}\sin{i}\\
-\sin{\zeta}\sin{i}&\cos{\zeta}\sin{i}&\cos{i}
\end{array}
\right)$$
$$\widetilde\RR_j=\left(
\begin{array}{ccc}
1-\sin^2(\zeta+\g)(1-\cos{i_j})&\sin{(\zeta+\g)}\cos{(\zeta+\g)}(1-\cos{i_j})&s_j\sin{(\zeta+\g)}\sin{i_j}\\
\sin{(\zeta+\g)}\cos{(\zeta+\g)}(1-\cos{i_j})&1-\cos^2(\zeta+\g)(1-\cos{i_j})&-s_j\cos{(\zeta+\g)}\sin{i_j}\\
-s_j\sin{(\zeta+\g)}\sin{i_j}&s_j\cos{(\zeta+\g)}\sin{i_j}&\cos{i_j}
\end{array}
\right)\ .$$}
Now, defining
\beqa{rps inversions}
\widetilde\cR_0= \widetilde\RR_0\circ\phi_{{\rm rps}^\complex}^{\rm jrd}\ ,\quad  \widetilde\cR_j=\widetilde\RR_j\circ\phi_{{\rm rps}^\complex}^{\rm jrd}\ ,\quad 
 \widetilde y^\ppj_{{\rm pl}}=\widetilde{\rm y}^\ppj_{\rm pl}\circ\phi_{{\rm rps}^\complex}^{\rm jrd}\ ,\quad 
 \widetilde y^\ppj_{{\rm pl}}=\widetilde{\rm y}^\ppj_{\rm pl}\circ\phi_{{\rm rps}^\complex}^{\rm jrd}
\eeqa
we have, from \equ{new formulae*},
\beqa{RPS}
 \arr{
x_{{\rm rps}^\complex}^\ppj=\widetilde\cR_0\widetilde\cR_j\widetilde x^\ppj_{\rm pl}\\
y_{{\rm rps}^\complex}^\ppj=\widetilde\cR_0\widetilde\cR_j\widetilde y^\ppj_{\rm pl}
}
\eeqa
We remark that the explicit expressions of $\widetilde\cR_0$, $\widetilde\cR_j$, $\widetilde x^\ppj_{\rm pl}$ and $\widetilde y^\ppj_{\rm pl}$ will be obtained from $\widetilde\RR_0$, $\widetilde\RR_j$, $\widetilde{\rm x}^\ppj_{\rm pl}$ and $\widetilde{\rm y}^\ppj_{\rm pl}$ replacing the appropriate arguments as
    \beqa{rps complex explicit}
&&\left\{
\begin{array}{lll}
\GG_j=\L_j-\ii \widetilde t_j\widetilde t_j^*\\
 i=\cos^{-1}\Big(\frac{\L_1+\L_2-\ii \widetilde t_1\widetilde t_1^*-\ii \widetilde t_2\widetilde t_2^*-\ii \widetilde t_3\widetilde t_3^*-\ii\widetilde T\widetilde T^*}{\L_1+\L_2-\ii \widetilde t_1\widetilde t_1^*-\ii \widetilde t_2\widetilde t_2^*-\ii \widetilde t_3\widetilde t_3^*}\Big)\\
1-\cos i_{1}=\frac{\big(\L_2-\ii \widetilde t_2\widetilde t_2^*\big)^2-
\big(\L_2-\ii \widetilde t_2\widetilde t_2^*-\ii \widetilde t_3\widetilde t_3^*\big)^2
}{2\big(\L_1+\L_2-\ii \widetilde t_1\widetilde t_1^*-\ii \widetilde t_2\widetilde t_2^*-\ii \widetilde t_3\widetilde t_3^*\big)\big(\L_1-\ii \widetilde t_1\widetilde t_1^*\big)}\\
\sin i_{1}=\frac{\sqrt{1-\Big(
(\L_1+\L_2-\ii  \widetilde t_1 \widetilde t_1^*-\ii  \widetilde t_2 \widetilde t_2^*-\ii  \widetilde t_3 \widetilde t_3^*)^2+(\L_1-\ii  \widetilde t_1 \widetilde t_1^*)^2-(\L_2-\ii  \widetilde t_2 \widetilde t_2^*)^2
\Big)^2}}{2\big(\L_1+\L_2-\ii  \widetilde t_1 \widetilde t_1^*-\ii  \widetilde t_2 \widetilde t_2^*-\ii  \widetilde t_3 \widetilde t_3^*\big)\big(\L_1-\ii  \widetilde t_1 \widetilde t_1^*\big)}\\
1-\cos i_{2}=\frac{\big(\L_1-\ii \widetilde t_1\widetilde t_1^*\big)^2-
\big(\L_1-\ii \widetilde t_1\widetilde t_1^*-\ii \widetilde t_3\widetilde t_3^*\big)^2
}{2\big(\L_1+\L_2-\ii \widetilde t_1\widetilde t_1^*-\ii \widetilde t_2\widetilde t_2^*-\ii \widetilde t_3\widetilde t_3^*\big)\big(\L_2-\ii \widetilde t_2\widetilde t_2^*\big)}\\
\sin i_{2}=\frac{\sqrt{1-\Big(
(\L_1+\L_2-\ii  \widetilde t_1 \widetilde t_1^*-\ii  \widetilde t_2 \widetilde t_2^*-\ii  \widetilde t_3 \widetilde t_3^*)^2+(\L_2-\ii  \widetilde t_2 \widetilde t_2^*)^2-(\L_1-\ii  \widetilde t_1 \widetilde t_1^*)^2
\Big)^2}}{2\big(\L_1+\L_2-\ii  \widetilde t_1 \widetilde t_1^*-\ii  \widetilde t_2 \widetilde t_2^*-\ii  \widetilde t_3 \widetilde t_3^*\big)\big(\L_2-\ii  \widetilde t_2 \widetilde t_2^*\big)}\\
\end{array}
\right.
 \nonumber\\
&& \left\{
\begin{array}{lll}
\zeta=\arg\frac{\widetilde T}{\sqrt{\ii \widetilde T \widetilde T^*}}\\
\g+\zeta=\arg\frac{\widetilde t_3}{\sqrt{\ii \widetilde t_3\widetilde t_3^*}}\\ 
\g_j+\g+\zeta=\arg\frac{\widetilde t_j}{\sqrt{\ii \widetilde t_j\widetilde t_j^*}}\\
 \ell_j=\l_j-\arg\frac{\widetilde t_j}{\sqrt{\ii \widetilde t_j\widetilde t_j^*}}
\end{array}
\right.
\eeqa
where we have used~\equ{D incli} to evaluate $\sin i_2$ and 
\beqa{i1}1-\cos i_1=\frac{\GG_2^2-(\GG-\GG_1)^2}{2\GG\GG_1}\ ,\qquad 1-\cos i_2=\frac{\GG_1^2-(\GG-\GG_2)^2}{2\GG\GG_2}\ .\eeqa
in turn implied by~\equ{D incli}.

\subsubsection{Analytical expression of $\phi_{{\rm rps}^\complex_\p}$}
To obtain the explicit expression of $\phi_{{\rm rps}_\p}$, we adopt a similar procedure as in the previous section. In this case, we rewrite~\equ{yx} in the form
\beqa{jrd++} 
\arr{
 x^\ppj_{\rm jrd}=\RR_0\RR_j\widetilde{\rm x}^\ppj_{\rm pl}\\
y^\ppj_{\rm jrd}=\RR_0\RR_j\widetilde{\rm y}^\ppj_{\rm pl}
  }\eeqa
 where
\beqano
&&\RR_0=\widetilde\RR_0\ ,\quad \RR_1=\widetilde\RR_1\ ,\quad  \widetilde{\rm x}^\ppu_{\rm pl}=\widetilde{\rm x}^\ppu_{\rm pl}\ ,\quad \widetilde{\rm y}^\ppu_{\rm pl}=\widetilde{\rm y}^\ppu_{\rm pl}\ ,\quad \RR_2=R^-_{313}(\zeta+\g, \s_2i_2)\nonumber\\
&&\widetilde{\rm x}^\ppd_{\rm pl}=R_3(\zeta+\g-\g_2)x^\ppd_{\rm orb}(-\L_2,\GG_2,-\ell_2)\ ,\quad  \widetilde{\rm y}^\ppd_{\rm pl}=R_3(\zeta+\g-\g_2)y^\ppd_{\rm orb}(-\L_2,\GG_2,-\ell_2)
\eeqano
 with $\widetilde\RR_0$, $\widetilde\RR_1$, $\widetilde{\rm x}_{\rm pl}^\ppu$, $\widetilde{\rm y}_{\rm pl}^\ppu$ as in~\equ{new formulae*}--\equ{rps inversions}, while, if $\P_2^-:=\left(
\begin{array}{lll}
1&0&0\\
0&-1&0\\
0&0&1
\end{array}
\right)$,
then
    \beqano
R_{313 }^-(\a,\b)&:=&R_3(\a)R_1(\b)R_3(\a)\P_2^-\nonumber\\
&=&\left(
\begin{array}{ccc}
1-\sin^2\a(1+\cos{\b})&\sin{\a}\cos{\a}(1+\cos{\b})&\sin{\a}\sin{\b}\\
\sin{\a}\cos{\a}(1+\cos{\b})&1-\cos^2\a(1+\cos{\b})&-\cos{\a}\sin{\b}\\
\sin{\a}\sin{\b}&-\cos{\a}\sin{\b}&\cos{\b}
\end{array}
\right)\ .
\eeqano
Note that here we have  used
    \beqano
&&\P_2^-R_3(g)x^\ppd_{\rm orb}(\L_2,\GG_2,\ell_2)=R_3(-g)x^\ppd_{\rm orb}(\L_2,\GG_2,-\ell_2)= x^\ppd_{\rm pl}(\L_2,\GG_2, -\ell_2,- g)=x^\ppd_{\rm pl}(-\L_2,\GG_2, -\ell_2,- g)\nonumber\\
&&\P_2^-R_3(g)y^\ppd_{\rm orb}(\L_2,\GG_2,\ell_2)=-R_3(-g)y^\ppd_{\rm orb}(\L_2,\GG_2,-\ell_2)= y^\ppd_{\rm pl}(-\L_2,\GG_2, -\ell_2,- g)\eeqano
Then we may write, for $\phi_{{\rm rps}^\complex_\p}$, the  expression
\beqa{RPSp***}
 \arr{
x_{{\rm rps}_\p^\complex}^\ppj=\cR_0\cR_j x^\ppj_{\rm pl}\\
y_{{\rm rps}_\p^\complex}^\ppj=\cR_0\cR_j y^\ppj_{\rm pl}
}
\eeqa
where
\beqano
\cR_0=\RR_0\circ\phi_{{{\rm rps^\complex_\p}}}^{\rm jrd}\ ,\quad \cR_j=\RR_j\circ\phi_{{{\rm rps^\complex_\p}}}^{\rm jrd}\ ,\quad x^\ppj_{\rm pl}=\widetilde{\rm x}^\ppj_{\rm pl}\circ\phi_{{{\rm rps^\complex_\p}}}^{\rm jrd}\ ,\quad y^\ppj_{\rm pl}=\widetilde{\rm y}^\ppj_{\rm pl}\circ\phi_{{{\rm rps^\complex_\p}}}^{\rm jrd}
\eeqano
The explicit form of $\phi_{{{\rm rps^\complex_\p}}}$ is obtained replacing the appropriate arguments as (as it follows from~\equ{PR})
\beqa{rps+ explicit}
&&\left\{
\begin{array}{lll}
\GG_1=\L_1-\ii  t_1 t_1^*\\
\GG_2=\L_2+\ii  t_2 t_2^*\\
i=\cos^{-1}\Big(\frac{\L_1-\L_2-\ii  t_1 t_1^*-\ii  t_2 t_2^*-\ii  t_3 t_3^*-\ii T T^*}{\L_1-\L_2-\ii  t_1 t_1^*-\ii  t_2 t_2^*-\ii  t_3 t_3^*-}\Big)\\
1-\cos i_{1}=\frac{\big(\L_2-\ii  t_2 t_2^*\big)^2-
\big(\L_2-\ii  t_2 t_2^*-\ii  t_3 t_3^*\big)^2
}{2\big(\L_1-\L_2-\ii  t_1 t_1^*-\ii  t_2 t_2^*-\ii  t_3 t_3^*\big)\big(\L_1-\ii  t_1 t_1^*\big)}\\
\sin i_{1}=\frac{\sqrt{1-\Big(
(\L_1-\L_2-\ii  t_1 t_1^*-\ii  t_2 t_2^*-\ii  t_3 t_3^*)^2+(\L_1-\ii  t_1 t_1^*)^2-(-\L_2-\ii  t_2 t_2^*)^2
\Big)^2}}{2\big(\L_1-\L_2-\ii  t_1 t_1^*-\ii  t_2 t_2^*-\ii  t_3 t_3^*\big)\big(\L_1-\ii  t_1 t_1^*\big)}\\
1+\cos i_{2}=\frac{\big(\L_1-\ii  t_1 t_1^*\big)^2-
\big(\L_1-\ii  t_1 t_1^*-\ii  t_3 t_3^*\big)^2
}{2\big(\L_1-\L_2-\ii  t_1 t_1^*-\ii  t_2 t_2^*-\ii  t_3 t_3^*\big)\big(-\L_2-\ii  t_2 t_2^*\big)}\\
-\sin i_{2}=\frac{\sqrt{1-\Big(
(\L_1-\L_2-\ii  t_1 t_1^*-\ii  t_2 t_2^*-\ii  t_3 t_3^*)^2+(-\L_2-\ii  t_2 t_2^*)^2-(\L_1-\ii  t_1 t_1^*)^2
\Big)^2}}{2\big(\L_1-\L_2-\ii  t_1 t_1^*-\ii  t_2 t_2^*-\ii  t_3 t_3^*\big)\big(-\L_2-\ii  t_2 t_2^*\big)}
\end{array}
\right.
\nonumber\\
&& \left\{
\begin{array}{lll}
\zeta=\arg\frac{T}{\sqrt{\ii T T^*}}\\
\g+\zeta=\arg\frac{ t_3}{\sqrt{\ii  t_3 t_3^*}}\\ 
\g_1+\g+\zeta=\arg\frac{ t_1}{\sqrt{\ii  t_1 t_1^*}}\\
\g+\zeta-\g_2=\arg\frac{ t_2}{\sqrt{\ii  t_2 t_2^*}}\\
 \ell_1=\l_1-\arg\frac{ t_1}{\sqrt{\ii  t_1 t_1^*}}\\
 - \ell_2=-\l_2-\arg\frac{ t_2}{\sqrt{\ii  t_2 t_2^*}}
\end{array}
\right.
\eeqa
having used the former equation in~\equ{i1} and
$$1+\cos i_2=\frac{\GG_1^2-(\GG+\GG_2)^2}{2\GG(-\GG_2)}$$
again following from~\equ{D incli}.\\
\subsubsection{Proof of Proposition~\ref{signs1}}

Since the third respective components of $\widetilde{\rm x}^\ppj_{\rm pl}$, $\widetilde{\rm y}^\ppj_{\rm pl}$ in \equ{RPS} and of ${\rm x}^\ppj_{\rm pl}$, ${\rm y}^\ppj_{\rm pl}$ in \equ{RPSp***}
vanish, we may replace the matrices $\widetilde\cR_j$ in  \equ{RPS} and
 $\cR_j$ in \equ{RPSp***} with the matrices that are obtained
 truncating at $0$ their respective third columns.
Collecting 
the formulae obtained after such truncation  with~\equ{rps complex explicit} and~\equ{rps+ explicit}, 
 the thesis follows.}
 
 \subsubsection{\label{formulae of reversed map}Analytical expression of $\phi_{{\rm rps}_\p}$}
Proposition~\ref{signs1}  has, as a consequence, that, formally, the formulae of the real ${\rm rps}_\p$ map in~\equ{rps+ map} can be derived from the ones of the real ${\rm rps}$ map  {given in  Ref.\cite[Equations (4.3)--(4.8) and Appendix A]{chierchiaPi11b}} (see also Ref.\cite[Section 4.3]{pinzari-th09}), for the case $n=2$, letting
 \beqa{real changes}
 &&\L_1=\L_1\ ,\quad
 \L_2\to-\L_2\ ,\quad
 \l_1=\l_1\ ,\quad
   \l_2=-\l_2\ ,\quad
    \eta_1=\eta_1\ ,\quad
      \xi_1=\xi_1\nonumber\\
      && \eta_2\to\ii\eta_2\ ,\quad
       \xi_2\to-\ii\xi_2\ ,\quad
       p_1\to\ii p\ ,\quad
         q_1\to-\ii q\ ,\quad
        p_2\to P\ ,\quad
        q_2\to Q\eeqa
  The result is as follows:
\beqa{rps++}\left\{
\begin{array}{ll}
x_{{\rm rps}_\p}^\ppj=\cR_0
\cR_j
 x^\ppj_{{\rm pl}}
\\ 
y_{{\rm rps}_\p}^\ppj=\cR_0
\cR_j
 y^\ppj_{{\rm pl}}
\end{array}
\right.
\qquad j=1,\ 2
\eeqa
where $\cR_0$, $\cR_j$, $x^\ppj_{{\rm pl}} $ and $y^\ppj_{{\rm pl}}$ are defined as follows.
Let  $s_1=-s_2=1$;
\beqa{c}{\rm c}_2^*&:=&\frac{1}{2(\L_1-\L_2)-\sum_{i=1}^2s_i(\eta_i^2+\xi_i^2)+p^2+q^2}\nonumber\\
 {\rm c}_1^*&:=&\frac{2\L_{2}-(\eta_2^2+\xi_2^2)-\frac{( p^2+ q^2)}{2}}{\Big(2(\L_1-\L_2)-\sum_{i=1}^2s_i(\eta_i^2+\xi_i^2)+p^2+q^2\Big)(2\L_1-(\eta_1^2+\xi_1^2))}
\nonumber\\
{\rm c}_2&:=&\frac{2\L_1-(\eta_1^2+\xi_1^2)+\frac{( p^2+ q^2)}{2}}{\Big(2(\L_1-\L_2)-\sum_{i=1}^2s_i(\eta_i^2+\xi_i^2)+p^2+q^2\Big)(2\L_2-(\eta_2^2+\xi_2^2))}\nonumber\\
 {\rm s}_i^*&:=&(-1)^i\sqrt{{\rm c}_i^*\Big(2-( p_i^2+ q_i^2){\rm c}_i^*\Big)}\ ,\quad   {\rm s}_2:=\sqrt{{\rm c}_2\Big(2-(p^2+q^2){\rm c}_2\Big)}\ .\eeqa
 Then
{\small \beqano
 \cR_0&:=&\left(
\begin{array}{ccc}
1-  Q^2{\rm c}^*_2&-  PQ{\rm c}^*_2&-  Q{\rm s}^*_2\\
-  PQ{\rm c}^*_2&1-  P^2{\rm c}^*_2&- P{\rm s}^*_2\\
  Q{\rm s}^*_2&  P{\rm s}^*_2&1-(  P^2+  Q^2){\rm c}^*_2
\end{array}
\right)\ ,\ 
\cR_1:=\left(
\begin{array}{ccc}
1-  q^2{\rm c}^*_1& p  q{\rm c}^*_1&-  q{\rm s}^*_1\\
 p  q{\rm c}^*_1&1-  p^2{\rm c}^*_1& p{\rm s}^*_1\\
  q{\rm s}^*_1&  -p{\rm s}^*_1&1-(  p^2+  q^2){\rm c}^*_1
\end{array}
\right) \nonumber\\
\cR_2&:=&\left(
\begin{array}{ccc}
1-  q^2{\rm c}_2&  p  q{\rm c}_2&-  q{\rm s}_2\\
 p  q{\rm c}_2&1-  p^2{\rm c}_2&  p{\rm s}_2\\
  q{\rm s}_2&  -p{\rm s}_2&1-(  p^2+  q^2){\rm c}_2
\end{array}
\right)
\eeqano}
while
$$
 (\L_j, \l_j,\eta_j, \xi_j)\to \Big( y^\ppj_{{\rm pl}} (\L_j, \l_j,\eta_j, \xi_j),  x^\ppj_{{\rm pl}} (\L_j, \l_j,\eta_j, \xi_j)\Big)\quad j=1,\ 2
$$
is the ``2--reversed planar Poincar\'e map'', given by
{\small\beqano
 \Big( y^\ppj_{{\rm pl}} (\L_j, \l_j,\eta_j, \xi_j),  x^\ppj_{{\rm pl}} (\L_j, \l_j,\eta_j, \xi_j)\Big)=\arr{
\Big(\widetilde y^\ppu_{{\rm pl}} (\L_1, \l_1,\eta_1, \xi_1),  \widetilde x^\ppd_{{\rm pl}} (\L_1, \l_1,\eta_1, \xi_1)\Big)\ \qquad\quad j=1\\
\\
\Big(\P_{2}^-\widetilde y^\ppd_{{\rm pl}} (\L_2, \l_2,\eta_2, \xi_2), \P_{2}^-\widetilde x^\ppd_{{\rm pl}} (\L_2, \l_2,\eta_2, \xi_2)\Big) \ j=2
}\eeqano}
with $\P_2^-=\left(
\begin{array}{ccc}
1&0\\
0&-1
\end{array}
\right)$; having used the relations
\beqa{planar poinc}
&&\widetilde x^\ppd_{{\rm pl}} (-\L_2, -\l_2,\ii\eta_2, -\ii\xi_2)=\P_{2}^-\widetilde x^\ppd_{{\rm pl}} (\L_2, \l_2,\eta_2, \xi_2)\nonumber\\
&&\widetilde y^\ppd_{{\rm pl}} (-\L_2, -\l_2,\ii\eta_2, -\ii\xi_2)=\P_{2}^-\widetilde y^\ppd_{{\rm pl}} (\L_2, \l_2,\eta_2, \xi_2)\ .
\eeqa

\subsubsection{\label{DAlembert rules}D'Alembert rules in the retrograde problem}
We  let
\beqa{RPS Ham}
\HH_{{\rm rps}^\complex_\p}(\L,\l,t,t^*):=\HH\circ \phi_{{\rm rps}^\complex_\p}(\L,\l,t,t^*)&=&\hh_{\rm k}(\L_1,\L_2)+\m f_{{\rm rps}^\complex_\p}(\L,\l,t,t^*)
\eeqa
the Hamiltonian~\equ{helio} expressed in the complex ${\rm rps}_\p$ coordinates.  
Proposition~\ref{signs1} implies that the the symmetries that have beed discussed in Ref.\cite[Section 6]{chierchiaPi11b} or Ref.\cite[Eqs. (3.19)--(3.24)]{chierchiaPi11c} have
{the following} counterpart {in the case of the retrograde problem}.
\paragraph{Invariance by rotations around $\CC$}
The following group of symmetries holds:
 \beqa{orc DAlembert}
f_{{\rm rps}^\complex_\p}(\L,\l, t, t^*)&=&f_{{\rm rps}^\complex_\p} \Big(\L,\ \l_1+g,\l_2-g, \cS_g ( t, t^*) \Big)\ ,\qquad \forall\ g\in \torus\ ,\eeqa
with
$\cS_{g}( t_1, t_2, t_3, t^*_1, t^*_2, t^*_3)=( t_1e^{\ii g}, t_2e^{\ii g}, t_3e^{\ii g}, t^*_1e^{-\ii g}, t^*_2e^{-\ii g}, t^*_3e^{-\ii g})$.
It is implied by 
the invariance of $f_{{\rm rps}^\complex_\p}$ under rotations around $\CC$ (corresponding to the Hamiltonian flow of $\GG$ in \equ{G}, which, in complex coordinates, is $\GG=\L_1-\L_2-{\rm i}t_1t_1^*-{\rm i}t_2t_2^*-{\rm i}t_3 t_3^*$), which, {in terms of ${\rm rps}^\complex_\p$, are}
\beqa{Rg}
&&\cR_g:\quad (\L,\l, t, t^*)\to \Big(\L,\ \l_1+g,\l_2-g, \cS_g ( t, t^*), T e^{\ii g},T^* e^{-\ii g}\Big)
\eeqa

\paragraph{Invariance by reflections} 
The following identities hold
\beqa{rot and reflect*}
&& f_{{\rm rps}^\complex_\p}(\L,\l, t, t^*)=f_{{\rm rps}_\p}(\L,\frac{\p}{2}-\l_1, -\frac{\p}{2}-\l_2,\cS_{{}_{1\leftrightarrow 2}} ( t, t^*))\nonumber\\
&&f_{{\rm rps}^\complex_\p}(\L,\l, t, t^*)=f_{{\rm rps}_\p}(\L,\l,\cS_3^-( t, t^*))\ .
\eeqa
with
\beqano
\arr{
\cS_{{}_{1\leftrightarrow 2}}( t_1, t_2, t_3, t^*_1, t^*_2, t^*_3):=(- t^*_1,- t^*_2, t_3,- t_1,- t_2, t^*_3) \\ 
\cS_3^-( t_1, t_2, t_3, t^*_1, t^*_2, t^*_3):=( t_1, t_2,- t_3, t^*_1, t^*_2,- t^*_3)\ .
}
\eeqano
Such indentities
 are implied by the invariance of $f_{{\rm rps}_\p}$ under the  tranformations
\beqa{orc DAlembert1}
&&\cR_3^-:\quad\quad (\L,\l, t, t^*)\to(\L,\l,\cS_3^-( t, t^*), -T,-T^*)\nonumber\\
&&\cR_{{}_{1\leftrightarrow 2}}:\quad(\L,\l, t, t^*)\to (\L,\frac{\p}{2}-\l_1, -\frac{\p}{2}-\l_2,\cS_{{}_{1\leftrightarrow 2}} ( t, t^*),T ,T^* )
\eeqa
in turn corresponding, in Cartesian coordinates, respectively, to
 \beqa{reflect1}
&&\begin{array}{llllllll}
\cR_3^-:\quad &x^\ppi\to\big(x^\ppi_1,\ x^\ppi_2,\ -x^\ppi_3\big)\ ,\quad &y^\ppi\to\big(y^\ppi_1,\ y^\ppi_2,\ -y^\ppi_3\big)\\
\cR_{{}_{1\leftrightarrow 2}}:\quad &x^\ppi\to\big(x^\ppi_2,\ x^\ppi_1,\ x^\ppi_3\big)\ ,\quad &y^\ppi\to\big(-y^\ppi_2,\ -y^\ppi_1,\ -y^\ppi_3\big)\ .
\end{array}\eeqa

\begin{acknowledgments}
I wish to thank the managing editor and, especially, the anonymous referees, whose careful remarks
 helped me  to improve in a substantial way  the presentation of the results. Thanks also  to A.~Celletti, 
A.~Giorgilli and V.~Kaloshin for their interest.
 \\
This research was  supported by H2020 Excellent Science (Grant 677793 StableChaoticPlanetM).\end{acknowledgments}


\begin{thebibliography}{99}
\bibitem{Note3bis} The only work known to the author where Deprit's coordinates are mentioned, up to~Ref.\cite{pinzari-th09}, is~Ref.\cite{ferrerO94}. However, such case is not really exhaustive, because it deals with the three--body problem, in which case Deprit's coordinates are the same as the classical Jacobi reduction of the nodes,~Ref.\cite{jacobi1842}.
\end{thebibliography}

\bibliographystyle{plain}

\end{document}